\documentclass{article}

\usepackage{amssymb}
\usepackage{amsfonts}
\usepackage[dvips]{graphicx}
\usepackage{mathrsfs}
\usepackage{amsmath}
\usepackage{latexsym}
\usepackage{caption}
\usepackage{amsthm}

\usepackage[usenames,dvipsnames]{pstricks}
\usepackage{epsfig}

\usepackage{graphicx}
\usepackage{placeins}
\usepackage{subfigure}

\newtheorem{theorem}{Theorem}
\newtheorem{lemma}[theorem]{Lemma}

 \newcommand{\nl}{\newline}
 \newcommand{\dist}{{\rm dist}}

\newcommand{\R}{{\mathbb{R}}}

 \newcommand{\diver}{{\rm div}}

 \newcommand{\half}{\frac{1}{2}}
 
 \newcommand{\darr}[4]{{\left\{\begin{array}{ll}
   {#1}&{#2}\\[0.2cm]
   {#3}&{#4}
 \end{array}\right.}}
  \newcommand{\darrsp}[4]{{\left\{\begin{array}{ll}
   {#1}&{#2}\\[0.4cm]
   {#3}&{#4}
 \end{array}\right.}}
 \newcommand{\tarr}[6]{{\left\{\begin{array}{lll}
   {#1}&{#2}\\
   {#3}&{#4}\\
   {#5}&{#6}
\end{array}\right.}}

\newcommand{\ia}{({\rm i})}
\newcommand{\ib}{({\rm ii})}
\newcommand{\ic}{({\rm iii})}

\def\Frac{\displaystyle\frac}

\title{On the Hardy constant of non-convex planar domains: \\
the case of the quadrilateral}
\date{}

\author{
G. Barbatis\footnote{Department of Mathematics,
 University of Athens,  15784 Athens, Greece}
 \and A. Tertikas
 \footnote{Department of Mathematics,
 University of Crete, 71409 Heraklion, Greece and  \nl
Institute of  Applied and Computational Mathematics,
FORTH, 71110 Heraklion, Greece}
}

\begin{document}

\date{}

\maketitle

\begin{abstract}
\noindent
The Hardy constant of a simply connected domain $\Omega\subset\R^2$ is the best constant for the inequality
\[
\int_{\Omega}|\nabla u|^2dx \geq c\int_{\Omega} \frac{u^2}{{\rm dist}(x,\partial\Omega)^2}\, dx \; , \;\;\quad  u\in C^{\infty}_c(\Omega).
\]
After the work of Ancona where the universal lower bound 1/16 was obtained, there has been a substantial interest on computing or estimating the Hardy constant of planar domains. In this work we determine the Hardy constant of an arbitrary quadrilateral in the plane. In particular we show that the Hardy constant is the same as that of a certain infinite sectorial region which has been studied by E.B. Davies.
\end{abstract}

\vspace{11pt}


\vspace{6pt}

\section{Introduction}

In the 1920's Hardy established the following inequality \cite{h}:
\begin{equation}
\int_0^{\infty} u'(t)^2dt \geq \frac{1}{4}\int_0^{\infty} \frac{u^2}{t^2}dt \; , \;\; \mbox{ for all }u\in C^{\infty}_c(0,\infty).
\label{or}
\end{equation}
The constant $1/4$ is the best possible, and equality is not attained for any non-zero function in the appropriate Sobolev space.

Inequality (\ref{or}) immediately implies the following inequality on $\R^N_+=\R^{N-1}\times (0,+\infty)$:
\begin{equation}
\int_{\R^N_+}|\nabla u|^2dx \geq \frac{1}{4}\int_{\R^N_+} \frac{u^2}{x_N^2}dx \; , \;\; \mbox{ for all }u\in C^{\infty}_c(\R^N_+),
\label{or1}
\end{equation}
where again the constant $1/4$ is the best possible.
The analogue of (\ref{or1}) for a domain $\Omega\subset\R^N$ is
\begin{equation}
\int_{\Omega}|\nabla u|^2dx \geq \frac{1}{4}\int_{\Omega} \frac{u^2}{d^2}\, dx \; , \;\; \mbox{ for all }u\in C^{\infty}_c(\Omega),
\label{or2}
\end{equation}
where $d=d(x)=\dist(x,\partial\Omega)$. However, (\ref{or2}) is not true without geometric assumptions on $\Omega$. The typical assumption made for the validity of (\ref{or2}) is that $\Omega$ is convex \cite{da2}. A weaker geometric assumption introduced in \cite{bft} is that $\Omega$ is weakly mean convex, that is
\begin{equation}
-\Delta d(x)\geq 0 \; , \quad \mbox{ in }\Omega,
\label{c}
\end{equation}
where $\Delta d$ is to be understood in the distributional sense. Condition (\ref{c}) is equivalent to convexity when $N=2$ but strictly weaker than convexity when $N\geq 3$
\cite{ak}.

In the last years there has been a lot of activity on Hardy inequality and improvements of it under the convexity or weak mean convexity assumption on $\Omega$;
see \cite{bm,bft,hhl,fmt}. If no geometric assumptions are imposed on $\Omega$, then one can still obtain inequalities of similar type. If for example $\Omega$ is bounded with $C^2$ boundary then one can still have inequality (\ref{or2}) for all $u\in C^{\infty}_c(\Omega_{\epsilon})$ where
$\Omega_{\epsilon}=\{x\in\Omega : d(x)<\epsilon\}$, provided $\epsilon>0$ is small enough \cite{fmt}. In the same spirit, under the same assumptions on $\Omega$ it was proved in \cite{bm} that there exists $\lambda\in\R$ such that
\begin{equation}
\int_{\Omega}|\nabla u|^2dx + \lambda \int_{\Omega}u^2 dx
\geq \frac{1}{4}\int_{\Omega} \frac{u^2}{d^2}\, dx \; , \;\; \mbox{ for all }u\in C^{\infty}_c(\Omega).
\label{or3}
\end{equation}

More generally, it is well known that for any bounded Lipschitz domain $\Omega\subset\R^N$ there exists $c>0$ such that
\begin{equation}
\label{hohi}
\int_{\Omega}|\nabla u|^2 dx\geq c\int_{\Omega}\frac{u^2}{d^2}dx \;\; , \qquad \mbox{ for all }u\in C^{\infty}_c(\Omega).
\end{equation}
Following \cite{da1} we call the best constant $c$ of inequality (\ref{hohi}) the Hardy constant of the domain $\Omega$.

In two space dimensions Ancona \cite{an} using Koebe's 1/4 theorem discovered the following remarkable result: for any simply connected domain $\Omega\subset\R^2$ there holds
\begin{equation}
\int_{\Omega}|\nabla u|^2dx
\geq \frac{1}{16}\int_{\Omega} \frac{u^2}{d^2}\, dx \; , \;\; \mbox{ for all }u\in C^{\infty}_c(\Omega).
\label{ancona}
\end{equation}
This result is typical of two space dimensions: Davies \cite{da1} has proved that no universal Hardy constant exists in dimension $N\geq 3$.

From now on we concentrate on two space dimensions. Two questions arise naturally, and have already been posed in the literature \cite{la,da1,da2,b,laso}:

\begin{description}

\item[ {\rm (1)} ] Given a simply connected domain $\Omega\subset\R^2$ find (or obtain information about) the Hardy constant of $\Omega$.

\item[ {\rm (2)} ] Find the best uniform Hardy constant valid for all simply connected domains $\Omega\subset\R^2$. Moreover, determine whether there
are extremal domains, that is domains $\Omega$ whose Hardy constant coincides with the best uniform Hardy constant.

\end{description}

Laptev and Sobolev \cite{laso} established a more refined version of Koebe's theorem and obtained a Hardy inequality which takes account of a quantitative measure of non-convexity. In particular they proved that if any $y\in\partial\Omega$ is the vertex of an
infinite sector $\Lambda$ of angle $\theta\in [\pi,2\pi]$ independent of $y$ such that $\Omega\subset \Lambda$, then the constant $1/16$ of (\ref{ancona}) can be replaced by $\pi^2/4\theta^2$.  The convex case corresponds to $\theta=\pi$, in which case the theorem recovers the $1/4$ in the case of convexity.
Analogous results were obtained recently in \cite{avla,ab}.

Davies \cite{da1} studied problem (1) in the case of an infinite sector of angle $\beta$. He used the symmetry of the domain to reduce the computation of the Hardy constant to the study of a certain ODE; see (\ref{ode}) below. In particular he established the following two results, which are also valid for the circular sector of angle $\beta$:

(a) The Hardy constant is $1/4$ for all angles $\beta\leq \beta_{cr}$, where $\beta_{cr} \cong 1.546\pi$.

(b) For $\beta_{cr} \leq \beta\leq 2\pi$ the Hardy constant strictly decreases with $\beta$ and in the limiting case $\beta=2\pi$ the Hardy constant is
$\cong 0.2054$.

Our aim in this work is to answer questions (1) and (2) in the particular case where $\Omega$ is a quadrilateral. Since the Hardy constant for any convex domain is $1/4$ we restrict our attention to non-convex quadrilaterals. In this case there is exactly one non-convex angle $\beta$,
$\pi <\beta <2\pi$. As we will see, this angle plays an important role and determines the Hardy constant. Our result reads as follows:
\newtheorem*{thm}{Theorem}
\begin{thm}
Let $\Omega$ be a non-convex quadrilateral with non-convex angle $\pi<\beta <2\pi$. Then
\begin{equation}
\int_{\Omega}|\nabla u|^2 dx \geq c_{\beta}\int_{\Omega} \frac{u^2}{d^2}dx  \; , \quad u\in C^{\infty}_c(\Omega),
\label{mix}
\end{equation}
where $c_{\beta}$ is the unique solution of the equation
\begin{equation}
\label{tans1}
\sqrt{c_{\beta}}\tan\big( \sqrt{c_{\beta}}(\frac{\beta-\pi}{2})\big) = 2 \bigg(
\frac{\Gamma(\frac{3+\sqrt{1-4c_{\beta}}}{4})}
{\Gamma(\frac{1+\sqrt{1-4c_{\beta}}}{4})}
\bigg)^2,
\end{equation}
when $\beta_{cr}\leq\beta < 2\pi$ and $c_{\beta}=1/4$ when $\pi <\beta\leq \beta_{cr}$. The constant $c_{\beta}$ is the best possible.
\end{thm}
As we shall see, the constant $c_{\beta}$ is precisely the Hardy constant of the sector of angle $\beta$, so equation (\ref{tans1}) provides an analytic description of the Hardy constant computed in \cite{da1} numerically. From (\ref{tans1}) we also deduce that the critical angle $\beta_{cr}$ in (b) is the unique solution in $(\pi,2\pi)$ of the equation
\begin{equation}
\label{tans2}
\tan\big( \frac{\beta_{cr}-\pi}{4}\big) = 4\bigg(\frac{\Gamma(\frac{3}{4})}{\Gamma(\frac{1}{4})}\bigg)^2.
\end{equation}
Relation (\ref{tans2}) was also obtained, amongst other interesting results, by Tidblom in \cite{ti}.
We also note that the constant $c_{2\pi}$ is the uniform Hardy constant for the class of all quadrilaterals.
The sharpness of the constant $c_{\beta}$ follows from the results of Davies \cite{da1}.

An important ingredient in the proof of our theorem
is the following elementary inequality valid on any domain $U$. Suppose $\partial U=\Gamma \cup\tilde\Gamma$.  Then, under certain assumptions, for any function
$\phi>0$ on $U\cup \Gamma$ we have
\begin{equation}
\int_U |\nabla u|^2dx    \geq -\int_U\frac{\Delta\phi}{\phi}u^2dx    +\int_{\Gamma} u^2 \frac{\nabla\phi}{\phi} \cdot\vec{\nu} dS
\label{lib1}
\end{equation}
for all smooth functions $u$ which vanish near $\tilde\Gamma$. Inequality (\ref{lib1}) will be applied to suitable subdomains $U_i$ of $\Omega$ and for
suitable choices of functions $\phi$.
Roughly, each subdomain $U_i$ consists of points whose nearest boundary point belongs to a different part of $\partial\Omega$.
The contribution along the boundary $\partial\Omega$ is zero because of the Dirichlet boundary conditions whereas there are non-zero interior boundary contributions that have to be taken into account.

The structure of the paper is simple: in Section \ref{sec:lemmas} we establish a number of auxiliary results that are used in Section \ref{sec:quads} where our theorem is proved.

\section{Auxiliary estimates}
\label{sec:lemmas}

Let $\beta >\pi$ be fixed. We start by defining the potential $V(\theta)$, $\theta\in (0,\beta)$,
\begin{equation}
V(\theta)=\tarr{\;\Frac{1}{\sin^2\theta} ,}{ 0< \theta  <\frac{\pi}{2} ,}{1,}{ \frac{\pi}{2}<\theta <\beta-\frac{\pi}{2},}
{\Frac{1}{\sin^2(\beta-\theta)},}{\beta-\frac{\pi}{2} <\theta <\beta.}
\label{v}
\end{equation}
For $c>0$ we consider the following boundary-value problem:
\begin{equation}
\darr{ -\psi''(\theta) =c V(\theta)\psi(\theta),}{ 0\leq\theta\leq\beta ,}{ \psi(0)=\psi(\beta)=0\, }{}
\label{ode}
\end{equation}
It was proved in \cite{da1} that the largest positive constant $c$ for which (\ref{ode}) has a positive solution coincides with Hardy constant of the sector of angle $\beta$. Due to the symmetry of the potential $V(\theta)$ this also coincides with
the largest constant $c$ for which the following boundary value problem has a solution:
\begin{equation}
\darr{ -\psi''(\theta) =c V(\theta)\psi(\theta),}{ 0\leq\theta\leq\beta/2,}{ \psi(0)=\psi'(\beta/2)=0\, .}{}
\label{ode1}
\end{equation}
Due to this symmetry, we shall identify the solutions of problems (\ref{ode}) and (\ref{ode1}).

The largest angle $\beta_{cr}$ for which the Hardy constant is $1/4$ for
$\beta\in [\pi,\beta_{cr}]$ was computed numerically in \cite{da1} and analytically in \cite{ti} where (\ref{tans2}) was established;
the approximate value is $\beta_{cr} \cong 1.546\pi$.

We first study the following algebraic equation
\begin{equation}
\sqrt{c}\tan\big( \sqrt{c}(\frac{\beta-\pi}{2})\big) = 2 \bigg(
\frac{\Gamma(\frac{3+\sqrt{1-4c}}{4})}
{\Gamma(\frac{1+\sqrt{1-4c}}{4})}
\bigg)^2.
\label{eq:alg}
\end{equation}
We note that choosing in (\ref{eq:alg}) $c=1/4$ we obtain $\beta_{cr}$ which is given by (\ref{tans2}).
\begin{lemma}
For any $\beta\geq \beta_{cr}$ there exists a unique $c=c_{\beta}$ satisfying (\ref{eq:alg}).  Moreover the function $\beta\mapsto c_{\beta}$ is smooth and strictly decreasing for $\beta\geq \beta_{cr}$. In particular we have
\[
c_{2\pi}<c_{\beta}<\frac{1}{4}\; \mbox{ for all }\; \beta_{cr} <\beta<2\pi.
\]
\end{lemma}
{\em Note.} From (\ref{eq:alg}) we  obtain the numerical estimate $c_{2\pi}\cong 0.20536$ of \cite{da1}.

{\em Proof.} Setting $x=\sqrt{1-4c}$ equation (\ref{eq:alg}) takes the equivalent form´
\[
G(x,\beta):= \frac{1}{2}(1-x^2)^{1/4} \tan^{1/2}\big( (1-x^2)^{1/2}\frac{\beta-\pi}{4}\big) -\frac{\Gamma(\frac{3+x}{4})}{\Gamma(\frac{1+x}{4})} =0,
\]
where we are interested in the range $0\leq x <1$ and $\beta$ is such that
\[
(1-x^2)^{1/2}\frac{\beta-\pi}{4} <\frac{\pi}{2}.
\]
For this range of $x$ and $\beta$ we can easily see that $G(x,\beta)$ is $C^{\infty}$. We will apply the Implicit Function Theorem. We first note that
$G(0,\beta_{cr})=0$. Moreover a simple but tedious computation gives
\begin{eqnarray*}
\frac{\partial G}{\partial x}(x,\beta)&=& -  \frac{x(\beta-\pi)}{16(1-x^2)^{1/4}}
\frac{1+  \tan^2\Big(     (1-x^2)^{1/2}\frac{\beta-\pi}{4}\Big)}{\tan^{1/2}\Big(     (1-x^2)^{1/2}\frac{\beta-\pi}{4}\Big)}  \\
&&  -\frac{x}{4(1-x^2)^{3/4}} \tan^{1/2}\Big(     (1-x^2)^{1/2}\frac{\beta-\pi}{4}\Big)
 -\frac{\Gamma( \frac{3+x}{4})}{4\Gamma( \frac{1+x}{4})}  \bigg( \frac{\Gamma'(\frac{3+x}{4})}{\Gamma(\frac{3+x}{4})}  
- \frac{\Gamma'(\frac{1+x}{4})}{\Gamma(\frac{1+x}{4})}    \bigg).
\end{eqnarray*}
Since
\[
\frac{d}{dx}\big(  \frac{\Gamma'(x)}{\Gamma(x)}\big) =\sum_{n=0}^{\infty}\frac{1}{(x+n)^2} >0,
\]
we conclude that $\partial G/\partial x <0$ for all $(x,\beta)$ with $0\leq x <1$ and
\[
\beta_{cr} \leq \beta <\frac{2\pi}{\sqrt{1-x^2}} +\pi .
\]
We also easily see that $\partial G/\partial \beta >0$ in the above range of $x$, $\beta$. This implies the existence and uniqueness locally near $\beta=\beta_{cr}$. A standard argument then gives the global existence of a smooth, strictly increasing function $x=x(\beta)$ for $\beta\geq \beta_{cr}$.
The proof is concluding by substituting $c=\frac{1-x^2}{4}$. $\hfill\Box$

We next study the boundary value problem (\ref{ode1}). The solution will be expressed using the hypergeometric function
\[
F(a,b,c ; z):=\frac{\Gamma(c)}{\Gamma(a)\Gamma(b)}\sum_{n=0}^{\infty}\frac{\Gamma(a+n)\Gamma(b+n)}{\Gamma(c+n)}\frac{z^n}{n!}.
\]
\begin{lemma}
Let $\beta>\beta_{cr}$. The boundary value problem (\ref{ode1}) has a positive solution if and only if $c$ solves (\ref{eq:alg}). In this case the solution is given by
\[
\psi(\theta)=\darr{ \Frac{\sqrt{2}\cos\big(\sqrt{c}(\beta-\pi)/2\big)\sin^{\alpha}(\theta/2)\cos^{1-\alpha}(\theta/2) }{F(\half,\half, \alpha +\half ; \half)}F(\half,\half,\alpha +\half ; \sin^2(\frac{\theta}{2})),}
{0<\theta\leq \frac{\pi}{2},}{ \cos\big( \sqrt{c}(\frac{\beta}{2}-\theta) \big),}{ \frac{\pi}{2} <\theta \leq \frac{\beta}{2},}
\]
where $\alpha$ is the largest solution of $\alpha(1-\alpha)=c$. Moreover $\psi\in H^1_0(0,\beta)$.
\label{lem:largeb}
\end{lemma}
{\em Proof.} Clearly the function
\[
\psi(\theta)=\cos (\sqrt{c_{\beta}}(\frac{\beta}{2}-\theta))\, , \qquad \frac{\pi}{2}\leq \theta\leq \frac{\beta}{2} \, .
\]
is a positive solution of the differential equation in $(\pi/2 ,\beta/2)$ and satisfies the boundary condition $\psi'(\beta/2)=0$.
For $\theta\in (0,\pi/2)$ we set $\xi =\sin^2\theta/2$ and $y(\theta)=\sin^{\alpha}(\theta/2) \cos^{1-\alpha}(\theta/2) w(\xi)$ and we obtain after some computations that $w(\xi)$ solves the hypergeometric equation
\[
\xi(1-\xi)w_{\xi\xi} +(2\xi +\alpha -\frac{3}{2})w_{\xi} +\frac{1}{4}w=0 \; , \qquad 0<\xi<\frac{1}{2},
\]
the general solution of which is described via hypergeometric functions $F(\alpha,\beta,\gamma,\xi)$ and is well-defined for $|\xi|<1$; see \cite{pz,as} for details and various properties of the hypergeometric functions. We thus conclude that the general solution of the differential equation in (\ref{ode1}) is
\begin{eqnarray*}
 y(\theta) &=&c_1 \sin^{\alpha}(\frac{\theta}{2}) \cos^{1-\alpha}(\frac{\theta}{2})F(\frac{1}{2},\frac{1}{2},\alpha+\frac{1}{2};\sin^2(\frac{\theta}{2})) \\
&& +
c_2 \sin^{1-\alpha}(\frac{\theta}{2}) \cos^{1-\alpha}(\frac{\theta}{2})F(1-\alpha,1-\alpha,\frac{3}{2}-\alpha;\sin^2(\frac{\theta}{2})).
\end{eqnarray*}
In order to maximize $c$ we take $c_2=0$. The matching conditions at $\theta=\pi/2$ force $c$ to satisfy equation (\ref{eq:alg}) and determine $c_1$.
$\hfill\Box$

\begin{lemma}
Let $\pi <\beta \leq \beta_{cr}$. The largest value of $c$ so that the boundary value problem (\ref{ode1}) has a positive solution is $c=1/4$.
For $\beta=\beta_{cr}$ the solution is
\[
\psi(\theta)=\darrsp{ \Frac{\cos\big(\frac{\beta_{cr}-\pi}{4}\big)\sin^{1/2}\theta }{F(\half,\half, 1 ; \half)}F(\half,\half,  1; \sin^2(\frac{\theta}{2})),}
{0<\theta\leq \frac{\pi}{2},}{ \cos\big( \frac{1}{2}(\frac{\beta_{cr}}{2}-\theta) \big),}{ \frac{\pi}{2} <\theta \leq \frac{\beta_{cr}}{2}.}
\]
\label{lem:smallb}
\end{lemma}
{\em Proof.} Let $c=1/4$. Working as in the proof of Lemma \ref{lem:largeb} we find that the general solution of the differential equation (\ref{ode1}) in $(0,\pi/2)$ now is
\begin{eqnarray*}
y(\theta) &=&c_1 \sin^{1/2}(\frac{\theta}{2}) \cos^{1/2}(\frac{\theta}{2})F(\frac{1}{2},\frac{1}{2},1;\sin^2(\frac{\theta}{2})) \\
&& + c_2 \sin^{1/2}(\frac{\theta}{2}) \cos^{1/2}(\frac{\theta}{2})F(\frac{1}{2},\frac{1}{2},1;\sin^2(\frac{\theta}{2}))
\int_{\sin^2(\theta/2)}^{1/2}\frac{dt}{t(1-t)F^2(\frac{1}{2},\frac{1}{2},1;t)}.
\end{eqnarray*}
The matching conditions at $\theta=\pi/2$ determine $c_1$ and $c_2$. In order for $\psi$ to be positive it is necessary that $c_2\geq 0$. This turns out to be equivalent to
\[
4\frac{ \Gamma^2(\frac{3}{4})}{ \Gamma^2(\frac{1}{4})} \geq \tan(\frac{\beta -\pi}{4}).
\]
This implies that $\beta\leq \beta_{cr}$ and in the case  $\beta=\beta_{cr}$ we have $c_2=0$. $\hfill\Box$

For our purposes it is useful to write the solution of (\ref{ode1}) in case $\beta\geq\beta_{cr}$
as a power series
\begin{equation}
\label{ps}
\psi(\theta) =\theta^{\alpha}\sum_{n=0}^{\infty}a_n\theta^{n} \; ,
\end{equation}
where $\alpha$ is the largest solution of the equation $\alpha(1-\alpha)=c$ in case $\beta>\beta_{cr}$ and $\alpha=1/2$ when $\beta=\beta_{cr}$.
We normalize the power series setting $a_0=1$; simple computations then give
\begin{equation}
a_1=0 \;\; , \qquad  a_2 =-\frac{\alpha(1-\alpha)}{6(1+2\alpha)}  .
\label{asymptotics}
\end{equation}

For our analysis it will be important to study the following two auxiliary functions:
\begin{equation}
f(\theta) =\frac{\psi'(\theta)}{\psi(\theta)} \; , \qquad \theta\in (0, \beta) \; ,
\label{f}
\end{equation}
and
\begin{equation}
g(\theta) =\frac{\psi'(\theta)}{\psi(\theta)}\sin\theta \; , \qquad \theta\in (0, \beta) \; ,
\label{g}
\end{equation}
where $\psi$ is the normalized solution of (\ref{ode}) described in Lemmas \ref{lem:largeb} and \ref{lem:smallb}.
We note that these functions depend on $\beta$.
Simple computations show that they respectively solve the differential equations
\begin{equation}
f'(\theta) +f^2(\theta)+cV(\theta)=0 \; , \qquad 0<\theta <\beta
\label{def}
\end{equation}
and
\begin{equation}
\label{deg}
g'(\theta) =-\frac{1}{\sin\theta} \Big[ g(\theta)^2 -\cos\theta \, g(\theta) +c \Big] \; \; , \quad 0<\theta \leq\pi/2,
\end{equation}
where $c=c_{\beta}$.

\begin{lemma}
\label{lem:g} Let $\pi\leq\beta\leq 2\pi$. The function $g(\theta)$ is monotone decreasing on $(0,\pi/2]$.
\end{lemma}
{\em Proof.} In the case where $\pi\leq\beta\leq \beta_{cr}$ we have $c=1/4$ and therefore monotonicity follows at once from (\ref{deg}).
Suppose now that $\beta_{cr}\leq\beta\leq 2\pi$. Using the asymptotics (\ref{asymptotics}) we obtain
\begin{equation}
g(\theta)=\alpha + (2a_2 -\frac{\alpha}{6})\theta^2 +O(\theta^3) \;\; , \qquad \mbox{ as }\theta\to 0+.
\label{g:asymp}
\end{equation}
Now, by (\ref{deg}) $g(\theta)$ is monotone decreasing in $[\theta_0,\pi/2]$ where $\theta_0\in [0,\pi/2]$ is determined by $\cos^2\theta_0 =4c$. Let $\rho^{+}(\theta)$ denote the largest root of the equation $t^2 -\cos\theta \, t +c$, $0\leq\theta\leq\theta_0$.
We note that $g(0)=\rho^+(0)$, $g'(0)=0$ and (by (\ref{g:asymp})) $g''(0)<0$. Hence there exists an non-empty interval $(0,\theta^*)$ on which $g$ is strictly monotone decreasing and, therefore, $g(\theta)>\rho^+(\theta)$.
To prove that $g$ is monotone decreasing on the whole $[0,\pi/2]$, let us assume that it is not.
Then there exists a least positive $\theta_1$ such that $g'(\theta_1)=0$. We then have $g(\theta_1)=\rho^+(\theta_1)$.
But $(\rho^+)'<0$, hence $g(\theta)<\rho^+(\theta)$ for $\theta<\theta_1$ close enough to
$\theta_1$. This contradicts the definition of $\theta_1$. $\hfill\Box$

\begin{lemma}
Let $\pi\leq\beta\leq 2\pi$. For $\pi/2\leq\gamma\leq\pi$ let $\theta_1$ be the angle in $[0,\pi/2]$ determined by the relation
\begin{equation}
\label{thetasss}
\cot\theta_1 =\sin\gamma .
\end{equation}
Then there holds
\begin{equation}
\label{calll}
\frac{2+\cos\gamma}{1+\sin^2\gamma}f(\theta_1)\geq f(\frac{\pi}{2}) \; , \qquad \frac{\pi}{2}\leq\gamma\leq\pi.
\end{equation}
\label{lem:sss1}
\label{lem:last}
\end{lemma}
{\em Proof.} We define
\[
Q(\gamma)= \frac{2+\cos\gamma}{1+\sin^2\gamma}f(\theta_1).
\]
We will establish that $Q$ is a decreasing function in $[\pi/2,\pi]$.
An easy calculation gives
\[
Q'(\gamma) =\frac{\cos\gamma \; (2+\cos\gamma) }{  (1+\sin^2\gamma)^2}\bigg[ f(\theta_1)^2  -\frac{\sin\gamma (\cos^2\gamma +4\cos\gamma+2)}
{\cos\gamma (2+\cos\gamma)}f(\theta_1) +c(1+\sin^2\gamma) \bigg],
\]
where $\theta_1=\theta_1(\gamma)$, $\pi/2\leq\gamma\leq\pi$.

We first consider the interval where $-2+\sqrt{2}\leq\cos\gamma\leq 0$. For such $\gamma$ we have $\cos^2\gamma +4\cos\gamma+2 \geq 0$ and the result follows at once.

We next consider the case where $-1\leq\cos\gamma \leq -2+\sqrt{2}$. The discriminant $\Delta$ of the quadratic polynomial above is
\[
\Delta =\frac{ \sin^2\gamma(\cos^2\gamma +4\cos\gamma+2)^2 -4c\cos^2\gamma(1+\sin^2\gamma)(2+\cos\gamma)^2}{\cos^2\gamma(2+\cos\gamma)^2}.
\]
However, since
\[
\frac{d}{dt}( t^2-4t+2)^2 = 4( t^2-4t+2)(t-2)<0 \; , \qquad 2-\sqrt{2}\leq t\leq 1,
\]
we conclude that $( t^2-4t+2)^2\leq 1$ for $2-\sqrt{2}\leq t\leq 1$ and therefore
\[
\Delta \leq \frac{ (1-\cos^2\gamma) -4c\cos^2\gamma (2-\cos^2\gamma)(2+\cos\gamma)^2}{\cos^2\gamma(2+\cos\gamma)^2} , \qquad \mbox{ for }
-1\leq\cos\gamma \leq -2+\sqrt{2}.
\]
Next we shall prove that $(1-\cos^2\gamma) -4c\cos^2\gamma (2-\cos^2\gamma)(2+\cos\gamma)^2\leq 0$ for $-1\leq\cos\gamma \leq -2+\sqrt{2}$.
For this we set $t=-\cos\gamma$ and we define $w(t)= 1-t^2 -4ct^2(2-t^2)(2-t)^2$, $t>0$. We have
\[
w'(t)=-2t\bigg( 1+4c [-3t^4+10t^3-4t^2-12t+8] \bigg).
\]
Now, the function $p(t)=-3t^4+10t^3-4t^2-12t+8$ has derivative
\[
p'(t)=(t-1)(-12t^2 +18t +10) -2 \leq 0 \; , \qquad 0\leq t\leq 1.
\]
Therefore $1+4cp(t) \geq 1+4cp(1) =1-4c\geq 0$ for $0\leq t\leq 1$. This in turn implies that $w(t)$ decreases in $[0,1]$. But
\[
w(2-\sqrt{2})=4\sqrt{2}-5 -64c(5\sqrt{2}-7) <0,
\]
since $c> (4\sqrt{2}-5)/(64(5\sqrt{2} -7))\approx 0.1444$.
We thus  conclude that $w(t)\leq 0$ for
$2-\sqrt{2}\leq t\leq 1$, which in turn implies that $\Delta\leq 0$ for $-1\leq\cos\gamma \leq -2+\sqrt{2}$. Therefore $Q(\gamma)$ is decreasing also in this this interval. Since $Q(\pi)=f(\pi/2)$, the proof is complete. $\hfill\Box$

\begin{lemma}
Let $\pi\leq\beta\leq 2\pi$ and $\pi/2\leq\gamma\leq\pi$. For $\theta\in [\pi/2 , (3\pi/2)-\gamma]$ denote by $\theta_1=\theta_1(\theta)$ be the angle in $[0,\pi/2]$ uniquely determined by the relation
\begin{equation}
\label{thetas}
\cot\theta_1 =-\cos(\theta+\gamma).
\end{equation}
Then there holds
\begin{equation}
\label{cal}
f(\theta_1) \geq f(\theta) \frac{1+\cos^2(\theta+\gamma)}{2+\sin(\theta+\gamma)} \, ,
\qquad \frac{\pi}{2}\leq \theta \leq\frac{3\pi}{2} -\gamma \, .
\end{equation}
\label{lem:sss}
\end{lemma}
{\em Proof.} For $\theta =\pi/2$ the corresponding value $\theta_*=\theta_1(\pi/2)$ is the one given by (\ref{thetasss}) hence the result is a consequence of Lemma \ref{lem:last}.

To prove (\ref{cal}) we shall consider $\theta_1$ as the free variable so that $\theta=\theta(\theta_1)$ is given by (\ref{thetas}). Since $f(\theta_1)$ satisfies $f'(\theta_1) +f^2(\theta_1) +c/\sin^2\theta_1=0$, it suffices to show that the function
\[
h(\theta_1) := f(\theta) \frac{1+\cos^2(\theta+\gamma)}{2+\sin(\theta+\gamma)}  \;\qquad (\theta=\theta(\theta_1))
\]
satisfies
\begin{equation}
H(\theta_1):= h'(\theta_1) +h^2(\theta_1) +\frac{c}{\sin^2\theta_1}  \leq  0\;\; , \qquad \theta_* \leq\theta_1\leq\frac{\pi}{2}\; ,
\label{nicesong}
\end{equation}
where $\theta_*\in (0,\pi/2)$ is determined by $\cot\theta_* = \sin \gamma$.

We express $H(\theta_1)$ in terms of $f(\theta)$ and $f'(\theta)$; we also use the fact that, by (\ref{thetas}),
\[
\frac{d\theta_1}{d\theta} =-\frac{ \sin(\theta+\gamma)}{1+\cos^2(\theta+\gamma)}.
\]
Using (\ref{def}) and setting $\omega=\theta+\gamma$ we obtain after some simple computations that
\begin{eqnarray}
&& H(\theta_1)=\frac{1+\cos^2\omega}{\sin\omega (2+\sin\omega)^2} \bigg[ 2(1+\cos^2\omega)(1+\sin\omega)f^2(\theta) + \label{eis}\\
&&\qquad\qquad + \cos\omega(\sin^2\omega +4\sin\omega +2)f(\theta) +2c(1+\sin\omega)(2+\sin\omega)\bigg]. \nonumber
\end{eqnarray}
In brackets we have a quadratic polynomial of $f(\theta)$ whose discriminant is itself a polynomial $P(t)$ of $t=-\sin\omega \in [-\cos \gamma,1]\subseteq [0,1]$,
\[
P(t)=(1-t) \Big[ t^5 +(16c-7)t^4 +12(1-4c)t^3 +4t^2 +12(8c-1)t +4(1-16c)\Big]=:(1-t)Q(t)\, .
\]
We observe that $Q(0)<0$ and $Q(1)=2>0$; moreover
\begin{equation}
Q'(t) =5t^4 +4(16c-7)t^3+36(1-4c)t^2 +8t +12(8c-1).
\label{q'}
\end{equation}
Recall that $1/8 <c \leq 1/4$, hence all the summands in (\ref{q'}) are non-negative in $[0,1]$ with the exception of $4(16c-7)$.
Since $|4(16c-7)|=28-64c < 36(1-4c) + 8 + 12(8c-1)$, we conclude that $Q'>0$ in $[0,1]$.

The above considerations imply that there exists a unique $t_0\in (0,1)$ such that $P(t)<0$ in $(0,t_0)$ and $P(t)>0$ in $(t_0,1)$.
This immediately implies that $H(\theta_1)\leq 0$ in the range $0<t<t_0$.

For $t_0<t<1$ the quadratic polynomial in (\ref{eis}) has two roots of the same sign as the sign of $t^2-4t+2$. The equation $t^2-4t+2=0$ has solutions $2\pm\sqrt{2}$. It follows that the quadratic polynomial above has negative two roots when $\max\{ t_0,2-\sqrt{2}\}<t<1$. Since $f(\theta)>0$, $0<\theta<\beta/2$, we conclude once again that $H(\theta_1)\leq 0$ in this case as well.
But we easily check that $Q(2-\sqrt{2})<0$, which implies that
$\max\{ t_0,2-\sqrt{2}\}=t_0$. This completes the proof. $\hfill\Box$

\begin{lemma}
Let $\pi\leq \beta\leq 2\pi$. The following inequalities hold:
\begin{eqnarray*}
\ia && \mbox{If $0\leq\omega\leq\pi/4$ then}\\
&& \hspace{1.5cm} f(\theta)\sin\theta\cos(\theta +\omega) + \alpha\cos\omega \geq 0 \; , \quad 0\leq\theta\leq\frac{\pi}{2}, \\[0.2cm]
\ib&& \mbox{If $3\pi/2 -\beta\leq \omega\leq 2\pi-\beta$ then}\\
&& \hspace{1.5cm} f(\theta)\cos(\theta+\omega)+\alpha[1+\sin(\theta+\omega)]\geq 0 \; ,\quad
\frac{\pi}{2}\leq\theta\leq\beta -\frac{\pi}{2}    \, ,\\[0.2cm]
\ic && \mbox{If $0\leq \omega\leq 2\pi-\beta$ then}\\
&& \hspace{1.5cm} -f(\theta)\cos(\theta+\omega) +\alpha[1-\sin(\theta+\omega)] \geq 0 \; ,
\quad\frac{\pi}{2}\leq\theta \leq\beta-\frac{\pi}{2}.
\end{eqnarray*}
\label{lem:quad}
\end{lemma}
{\em Proof.} (i) The inequality is trivially true for $0\leq\theta\leq\pi/2 -\omega$, so we restrict our attention to the interval $\pi/2 -\omega\leq\theta\leq\pi/2$. We must prove that
\begin{equation}
f(\theta)\leq F(\theta) \; , \qquad \frac{\pi}{2} -\omega\leq\theta\leq\frac{\pi}{2},
\label{fg}
\end{equation}
where $f$ is given by (\ref{f}) and
\[
F(\theta)=-\alpha\frac{\cos\omega}{\sin\theta\cos(\theta+\omega)}.
\]
Using the fact that $\sqrt{c}\leq\alpha$ we have
\begin{eqnarray}
F(\frac{\pi}{2})-f(\frac{\pi}{2}) &=& \alpha \cot\omega-\sqrt{c}\tan[\sqrt{c}(\frac{\beta}{2}-\frac{\pi}{2})] \nonumber\\
&\geq& \alpha\Big\{ \cot\omega-\tan[\sqrt{c}(\frac{\beta}{2}-\frac{\pi}{2})]\Big\} \nonumber\\
&= & \frac{\alpha}{\sin\omega \cos[\sqrt{c}(\frac{\beta}{2}-\frac{\pi}{2})]}\cos\Big(\sqrt{c}(\frac{\beta}{2}-\frac{\pi}{2}) +\omega \Big) \nonumber\\
&\geq& 0,
\label{p2}
\end{eqnarray}
since $0<\sqrt{c}(\frac{\beta}{2}-\frac{\pi}{2}) +\omega \leq \frac{\beta}{4}-\frac{\pi}{4}+\omega\leq \pi/2$.

We shall prove that $F'(\theta) +F(\theta)^2 +\frac{c}{\sin^2\theta} \leq 0$ in
$[\pi/2-\omega,\pi/2]$.
This, combined with (\ref{def}) and (\ref{p2}) will imply that $f(\theta)\leq F(\theta)$ in $[\pi/2-\omega,\pi/2]$.

Recalling that $c= \alpha(1-\alpha)$, we have for $\theta\in [\pi/2-\omega,\pi/2]$,
\begin{eqnarray*}
&&F'(\theta)+F^2(\theta) +\frac{c}{\sin^2\theta} \\
&=&\frac{\alpha\cos\omega\cos\theta}{\sin^2\theta\cos(\theta+\omega)} -
\frac{\alpha\cos\omega\sin(\theta+\omega)}{\sin\theta\cos^2(\theta+\omega)}
 +\frac{\alpha^2\cos^2\omega}{\sin^2\theta\cos^2(\theta+\omega)} +\frac{c}{\sin^2\theta} \\
&=&\alpha\frac{\cos\omega\cos\theta\cos(\theta+\omega) -\cos\omega\sin(\theta+\omega)\sin\theta +\alpha\cos^2\omega +(1-\alpha)\cos^2(\theta+\omega)}{\sin^2\theta \cos^2(\theta+\omega)} \\
&=&\alpha\frac{ 2\cos\omega\cos\theta\cos(\theta+\omega) -(1-\alpha) [\cos^2\omega-\cos^2(\theta+\omega)]}{\sin^2\theta \cos^2(\theta+\omega)} \\
&=&\alpha\frac{ 2\cos\omega\cos\theta\cos(\theta+\omega) -(1-\alpha)\sin\theta \sin(\theta+2\omega)}{\sin^2\theta \cos^2(\theta+\omega)} \\
&\leq& 0,
\end{eqnarray*}
since the last term is the sum of two non-positive terms. Hence $\ia$ has been proved.

(ii) We first note that
\[
f(\theta)=\sqrt{c}\tan\Big( \sqrt{c}(\frac{\beta}{2}-\theta)\Big), \qquad\quad \frac{\pi}{2}\leq\theta\leq\beta -\frac{\pi}{2},
\]
and
\[
-\frac{\pi}{4}\leq \sqrt{c}(\frac{\beta}{2}-\theta) \leq\frac{\pi}{4} \;, \qquad \quad \frac{\pi}{2}\leq\theta\leq\beta -\frac{\pi}{2}.
\]
It follows that the required inequality is written equivalently,
\[
\alpha(1+\sin(\omega+\theta))\cos(\sqrt{c}( \frac{\beta}{2}-\theta)) +\sqrt{c}\sin(\sqrt{c}( \frac{\beta}{2}-\theta))\cos(\omega+\theta)\geq 0 \; ,
\;\; \frac{\pi}{2}\leq\theta\leq\beta -\frac{\pi}{2}.
\]
Hence, since $\alpha\geq\sqrt{c}$,
\begin{eqnarray}
&&\hspace{-2cm}\alpha(1+\sin(\theta+\omega))\cos(\sqrt{c}( \frac{\beta}{2}-\theta)) +\sqrt{c}\sin(\sqrt{c}( \frac{\beta}{2}-\theta))\cos(\theta+\omega) \nonumber\\
&\geq&\sqrt{c}\Big\{(1+\sin(\theta+\omega))\cos(\sqrt{c}( \frac{\beta}{2}-\theta)) +\sin(\sqrt{c}( \frac{\beta}{2}-\theta))\cos(\theta+\omega)\Big\} \nonumber\\
&=&\sqrt{c}\Big\{ \cos(\sqrt{c}( \frac{\beta}{2}-\theta)) +\sin[\sqrt{c}( \frac{\beta}{2}-\theta) +\theta+\omega]\Big\} \nonumber\\
&=&\sqrt{c}\Big\{\cos(\sqrt{c}( \frac{\beta}{2}-\theta)) -\cos[\frac{\pi}{2}+\sqrt{c}( \frac{\beta}{2}-\theta) +\theta+\omega]\Big\}.\nonumber\\
&=& 2\sqrt{c} \sin\Big[ \sqrt{c}(\frac{\beta}{2}-\theta)+\frac{\pi}{4}+\frac{\theta}{2}+\frac{\omega}{2}\Big] \sin(\frac{\pi}{4} +\frac{\theta}{2}+\frac{\omega}{2}).
\label{s1}
\end{eqnarray}
But for the given range of $\omega$ and $\theta$ we have
\[
0\leq \frac{\pi}{4} +\frac{\theta}{2}+\frac{\omega}{2} \leq\pi \qquad
\mbox{ and }\qquad
0\leq \sqrt{c}(\frac{\beta}{2}-\theta)+\frac{\pi}{4}+\frac{\theta}{2}+\frac{\omega}{2} \leq \pi.
\]
Hence the last quantity in (\ref{s1}) is non-negative.

(iii) We have $\cos(\theta+\omega)\leq 0$ for $\frac{\pi}{2}\leq\theta\leq\beta-\frac{\pi}{2}$, therefore the inequality is trivial for $\theta\in [\pi/2,\beta/2]$ (since $f \geq 0$ there).
We now consider the complementary interval $\beta/2\leq\theta\leq\beta -\pi/2$. Arguing as in (\ref{s1}) above we see that it suffices to prove that
\[
-\sin(\sqrt{c}( \frac{\beta}{2}-\theta))\cos(\theta+\omega) +[1-\sin(\theta+\omega) ]\cos(\sqrt{c}( \frac{\beta}{2}-\theta))\geq 0,
\]
or equivalently,
\begin{equation}
\cos\big(\sqrt{c}(\theta- \frac{\beta}{2})\big) \geq \sin \big(\sqrt{c}( \frac{\beta}{2}-\theta)+\theta+\omega\big) \; , \quad
\frac{\beta}{2}\leq\theta\leq\beta -\frac{\pi}{2}.
\label{ph}
\end{equation}
We have
\[
\cos\big(\sqrt{c}(\theta- \frac{\beta}{2})\big) -\sin \big(\sqrt{c}( \frac{\beta}{2}-\theta)+\theta+\omega \big)=
-2\sin\big(\frac{\pi}{4} -\frac{\theta+\omega}{2} \big)\sin\big(\sqrt{c}(\frac{\beta}{2}-\theta) +\frac{\theta+\omega}{2} -\frac{\pi}{4}\big)
\]
Since $\beta+\omega \leq 2\pi$, we have
\[
0\leq \frac{\theta}{2}+\frac{\omega}{2} -\frac{\pi}{4} \leq\frac{\pi}{2}
\]
and
\[
0\leq \sqrt{c}(\frac{\beta}{2}-\theta) +\frac{\theta+\omega}{2} -\frac{\pi}{4}
\leq -\frac{\sqrt{c}(\beta-\pi)}{2} +\frac{\beta+\omega}{2} \leq\frac{\pi}{2},
\]
hence (\ref{ph}) is true. $\hfill\Box$




\section{Proof of the Theorem}
\label{sec:quads}

In this section we will give the proof of our Theorem.
We start with a lemma that plays fundamental role in our argument and will be used repeatedly. We do not try to obtain the most general statement and for simplicity we restrict ourselves to assumptions that are sufficient for our purposes.

Let $U$ be a domain and assume that $\partial U=\Gamma\cup\tilde\Gamma$ where $\Gamma$ is Lipschitz continuous. We denote by $\vec{\nu}$ the exterior unit normal on $\Gamma$.
\begin{lemma}
Let $\phi\in H^1_{\rm loc}(U)$ be a positive function such that $\nabla\phi /\phi \in L^2(U)$ and $\nabla\phi /\phi$ has an $L^1$ trace on $\Gamma$ in the sense that $v \nabla\phi/\phi$ has an $L^1$ trace on $\partial U$ for every  $v\in C^{\infty}(\overline{U})$ that vanishes near $\tilde\Gamma$.
Then
\begin{equation}
\int_U |\nabla u|^2dx\, dy    \geq -\int_U\frac{\Delta\phi}{\phi}u^2dx\, dy    +\int_{\Gamma}  \frac{\nabla\phi}{\phi} \cdot\vec{\nu} u^2 dS
\label{lib}
\end{equation}
for all smooth functions $u$ which vanish near $\tilde\Gamma$. Here $\Delta\phi$ is understood in the distributional sense. \nl
If in particular there exists $c\in\R$ such that
\begin{equation}
-\Delta\phi\geq\frac{c}{d^2}\phi \; ,
\label{1}
\end{equation}
in the weak sense in $U$, where $d=\dist(x,\tilde\Gamma)$, then
\begin{equation}
\int_U |\nabla u|^2dx\,  dy  \geq c\int_{U}\frac{u^2}{d^2}dx\, dy   +\int_{\Gamma} u^2 \frac{\nabla\phi}{\phi} \cdot\vec{\nu} dS
\label{2}
\end{equation}
for all functions $u\in C^{\infty}(\overline{U})$ that vanish near $\tilde\Gamma$.
\label{lem:1}
\end{lemma}
{\em Proof.} Let $u$ be a function in $C^{\infty}(\overline{U})$ that vanishes near $\tilde\Gamma$. We denote $\vec{T}=-\nabla\phi/\phi$. Then
\begin{eqnarray*}
\int_U u^2 \diver \vec{T}\, dx\, dy &=& -2\int_U u\nabla u\cdot \vec{T}\, dx\, dy +\int_{\Gamma}u^2 \vec{T}\cdot \vec{\nu} \, dS  \\
&\leq &\int_U  |\vec{T}|^2u^2 dx\, dy +\int_U |\nabla u|^2 dx\, dy + \int_{\Gamma}u^2 \vec{T}\cdot \vec{\nu} \, dS\, ,
\end{eqnarray*}
that is
\[
\int_U |\nabla u|^2 dx\, dy \geq \int_U (\diver \vec{T} -|\vec{T}|^2) u^2 dx\, dy -\int_{\Gamma}  \vec{T}\cdot \vec{\nu}u^2 dS\, .
\]
Using assumption (\ref{1}) we obtain (\ref{2}). $\hfill\Box$

Let us now consider a non-convex quadrilateral $\Omega$, with vertices $O$, $A$, $B$ and $C$  (as in the diagrams) and corresponding angles $\beta$, $\gamma$, $\delta$ and $\zeta$. We assume that the non-convex vertex is $O$ and, is located at the origin, and that the side $OA$ lies along the positive $x$-axis and has length one.

Our argument depends fundamentally on two geometric features of the quadrilateral $\Omega$. While in all cases the methodology remains the same, the technical details are different.
The first feature is whether or not one of the angles adjacent to the non-convex one is larger than $\pi/2$. The second one is related to the structure of the equidistance curve
\[
\Gamma=\{ P \in\Omega : \dist(P, OA \cup OC) =\dist(P, AB\cup BC)\}.
\]
Clearly the curve $\Gamma$ consists of line and parabola segments. Taking also account of symmetries, each non-convex quadrilateral $\Omega$ fits within one of the following five types, each one of which will be dealt with separately:

\begin{figure}[ht]
\begin{minipage}[b]{0.5\linewidth}
\centering
\includegraphics[scale=0.4]{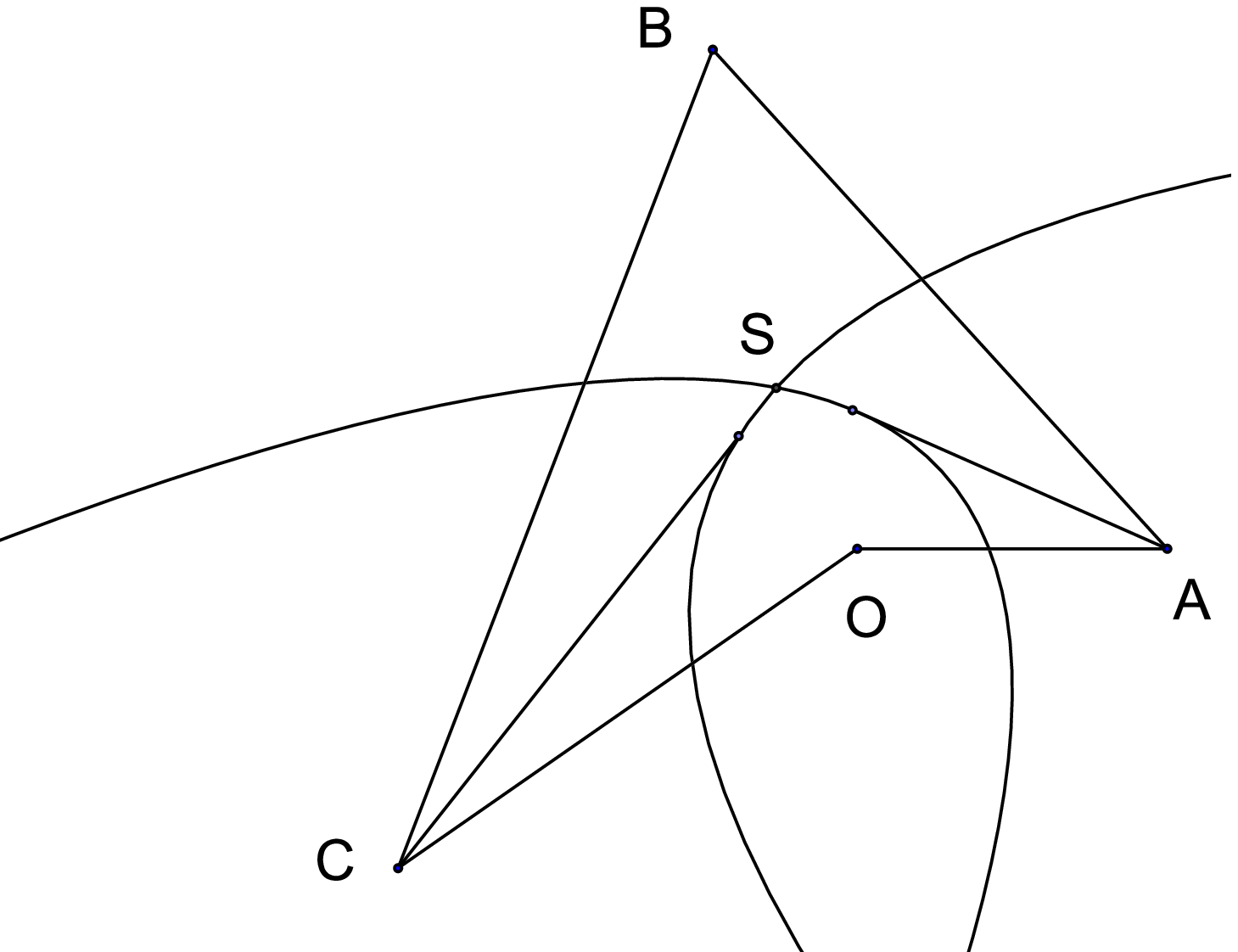}
\caption*{Type A1}
\end{minipage}
\hspace{-35pt}
\begin{minipage}[b]{0.5\linewidth}
\centering
\includegraphics[scale=0.4]{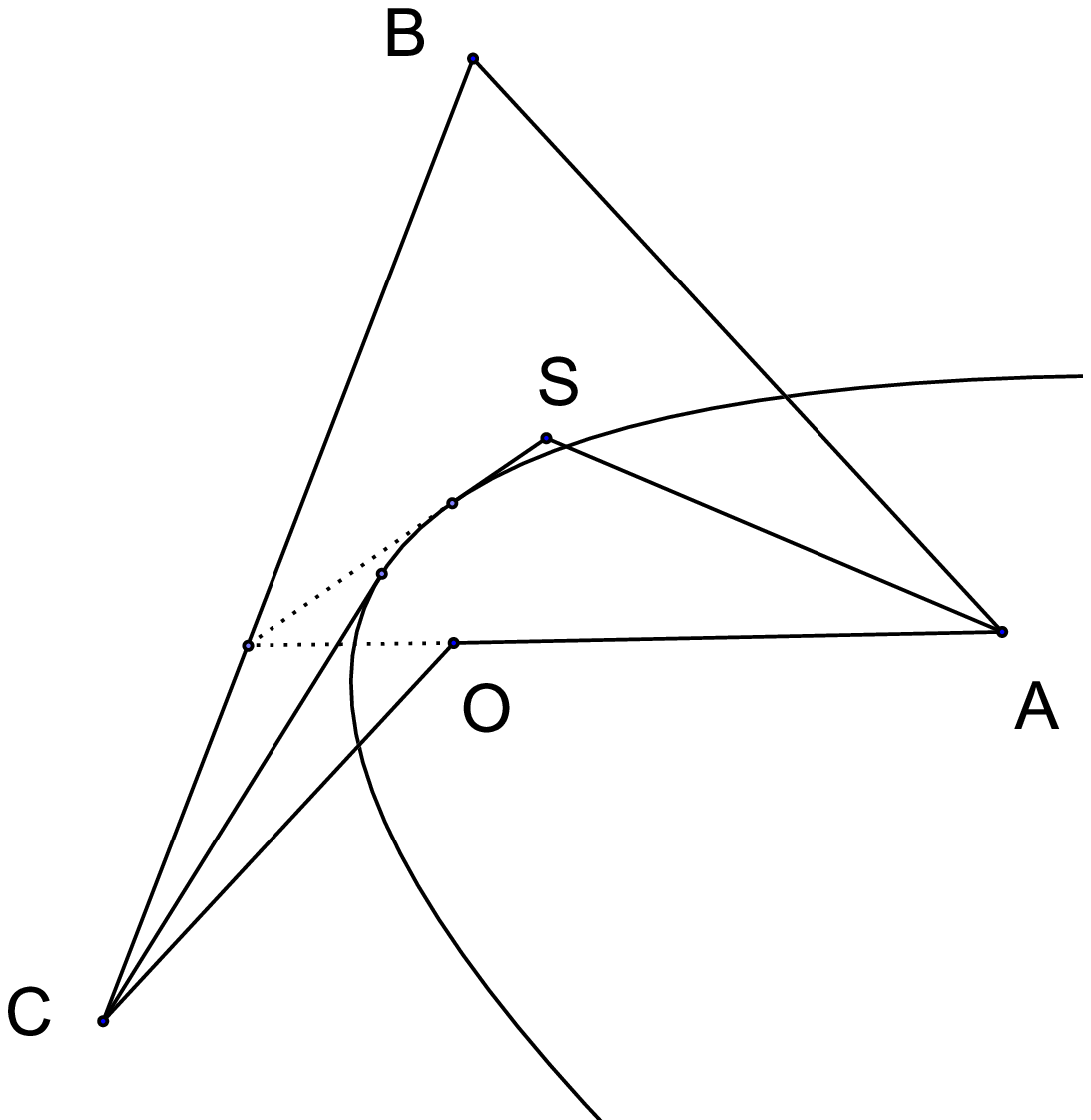}
\caption*{Type A2}
\end{minipage}
\end{figure}

{\bf Type A1.} We have $\gamma\leq\pi/2$, $\zeta\leq\pi/2$ and the curve $\Gamma$ consists of two line and two parabola segments
(Here we also include the special case where $\Gamma$ consists of two line segments and one parabola segment.)

{\bf Type A2.} We have $\gamma\leq\pi/2$, $\zeta\leq\pi/2$ and the curve $\Gamma$ consists of three line segments and one parabola segment.

\begin{figure}[ht]
\begin{minipage}[b]{0.5\linewidth}
\centering
\includegraphics[scale=0.4]{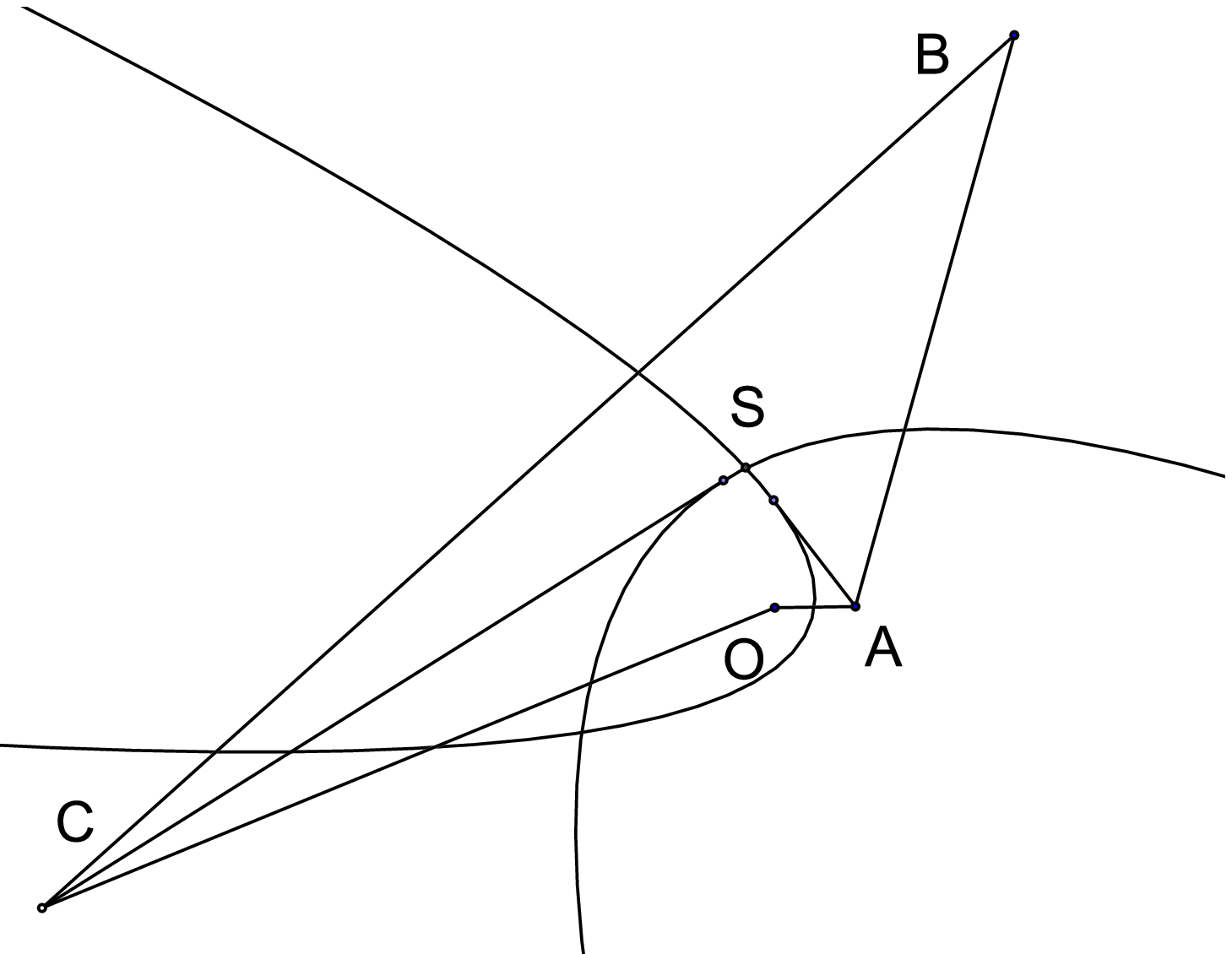}
\caption*{Type B1}
\end{minipage}
\hspace{-35pt}
\begin{minipage}[b]{0.5\linewidth}
\centering
\includegraphics[scale=0.4]{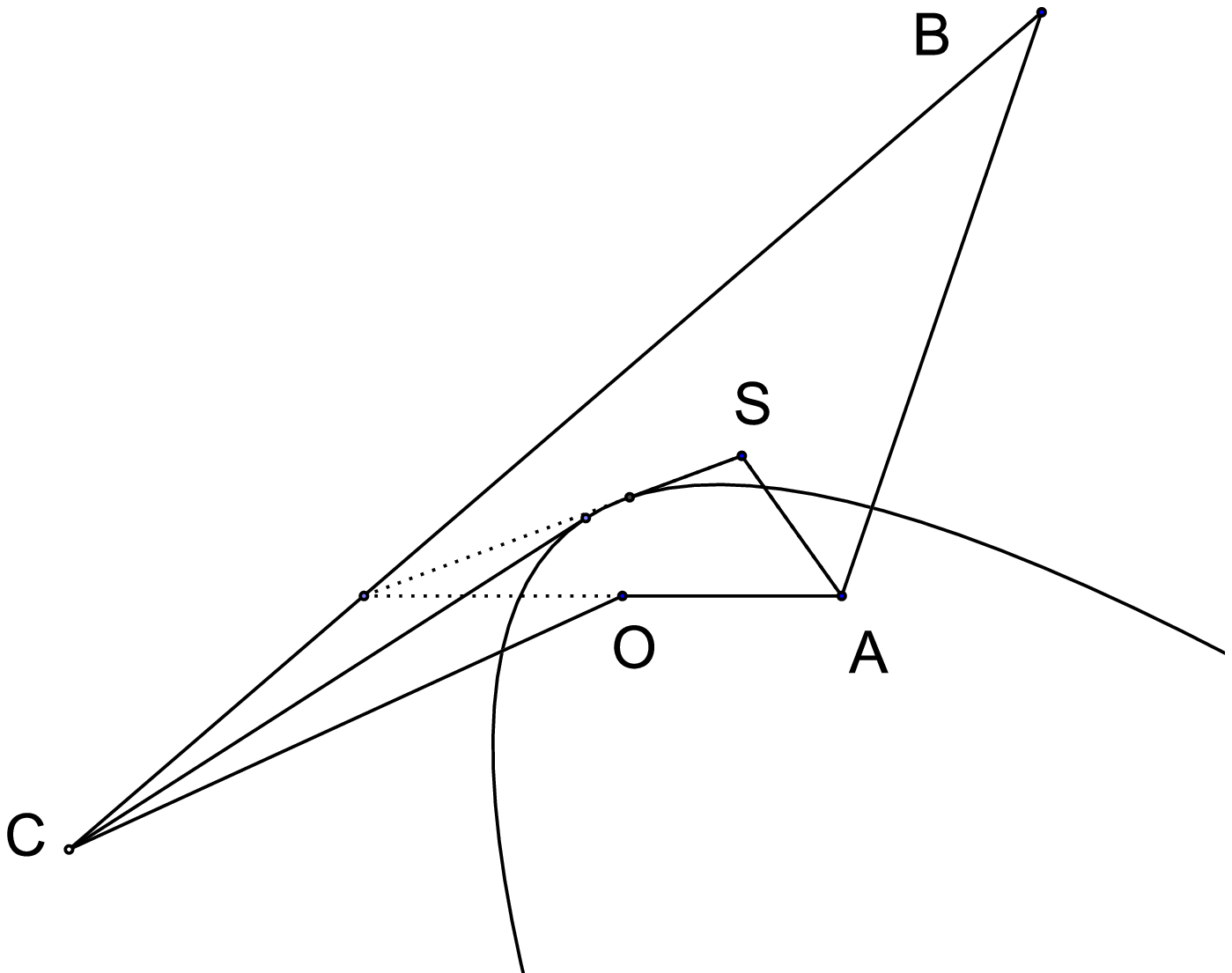}
\caption*{Type B2}
\end{minipage}
\end{figure}

\begin{center}
\begin{figure}
\begin{center}
\includegraphics[scale=0.32]{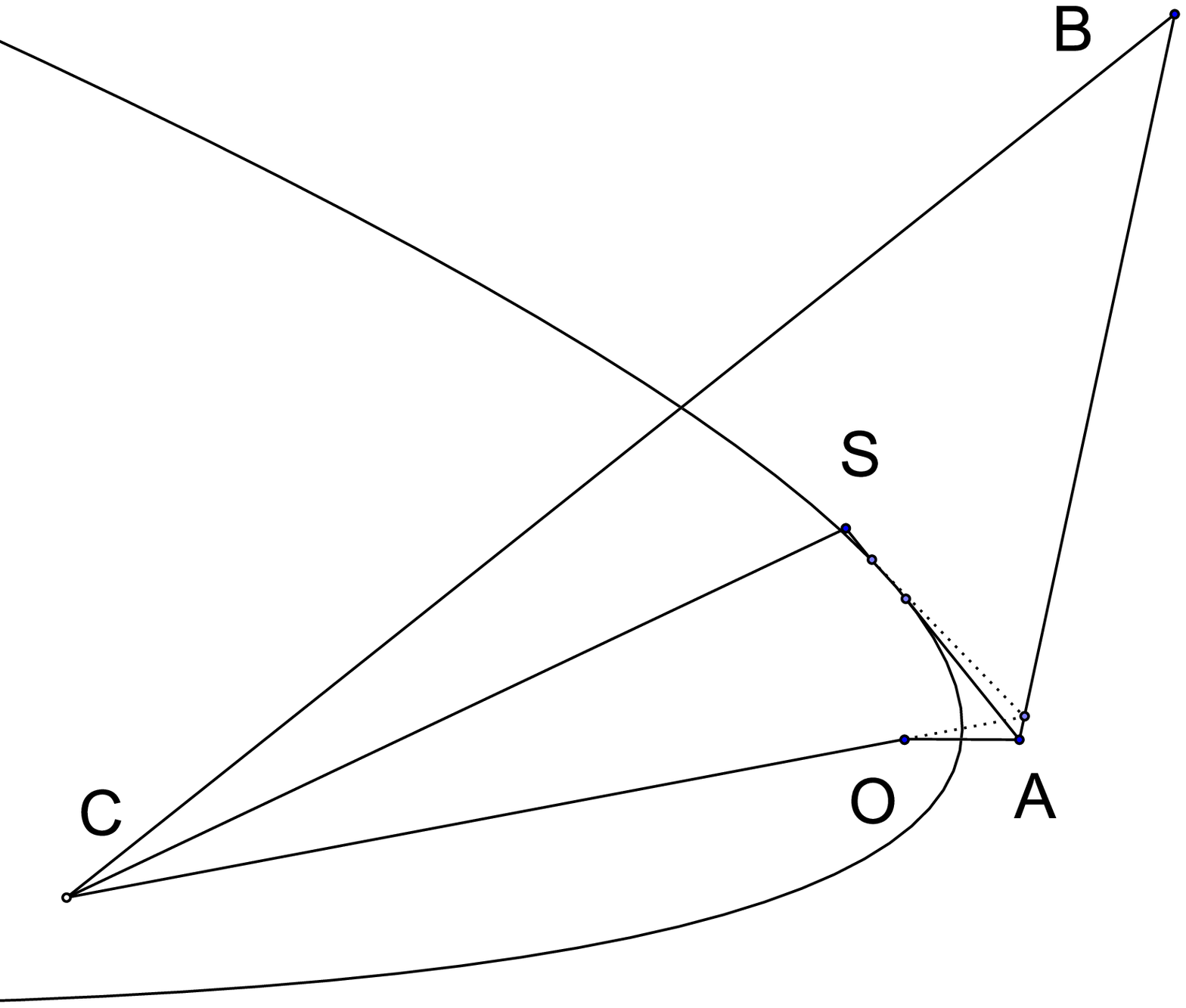}
\caption*{Type B3}
\end{center}
\end{figure}
\end{center}

{\bf Type B1.} $\gamma >\pi/2$ and the curve $\Gamma$ consists of two line segments and two parabola segments.
(Here we also include the special case where $\Gamma$ consists of two line segments segments and one parabola segment.)

{\bf Type B2.} $\gamma >\pi/2$ and the curve $\Gamma$ consists of three line and one parabola segment: starting from the point $A$ we first have two line segments, then a parabola segment and then a last line segment.

{\bf Type B3.} $\gamma >\pi/2$ and the curve $\Gamma$ consists again of three line and one parabola segment: starting from the point $A$ we first have a line segment, then a parabola segment and then two more line segments.

In all cases the curve $\Gamma$ divides $\Omega$ into two parts $\Omega^-$ and $\Omega^+$ where points in $\Omega^-$ have nearest boundary point on $OA \cup OC$ and points on $\Omega^+$ have nearest boundary points on $AB\cup BC$.
We denote by $\vec{\nu}$ the unit normal along $\Gamma$ which is outward with respect to $\Omega^-$.  We also denote by $S$ the point where $\Gamma$ intersects the bisector at the vertex $B$.

We shall often make use of the following simple fact: let $P$ be the parabola determined by the origin and the line $x\sin\alpha  +y\cos\alpha  +l=0$, where $l>0$. The exterior (with respect to the convex component)
unit normal along $\partial P$ is given in polar coordinates by
\begin{equation}
\vec{\nu}=\frac{(\cos\theta -\sin\alpha, \sin\theta-\cos\alpha)}{\sqrt{2-2\sin(\theta+\alpha)}}.
\label{normal}
\end{equation}

{\bf\em Proof of Theorem: type A1.} We parametrize $\Gamma$ by the polar angle $\theta\in [0,\beta]$. For $\theta\in [0,\pi/2]$ $\Gamma$ is a straight line;
the same is true for $\theta\in [\beta -\pi/2 ,\beta]$. Finally, for $\theta\in [\pi/2 ,\beta-\pi/2]$ $\Gamma$ consists of segments of two parabolas. These parabolas meet at the point $S$ which is equidistant from $AB$, $BC$ and the origin. Let $\theta_0$ be the polar angle of $S$.
We assume without loss of generality that $\theta_0\leq\beta /2$. Hence $\Gamma$ consists of four segments which when parametrized by the polar angle $\theta$ are described as
\[
\Gamma_1=\{ 0\leq \theta\leq \pi/2\} , \;
\Gamma_2=\{ \pi/2 \leq \theta\leq \theta_0 \} , \;
\Gamma_3=\{ \theta_0\leq \theta\leq \beta-\frac{\pi}{2}\} , \;
\Gamma_4=\{\beta-\frac{\pi}{2}\leq\theta \leq\beta\}.
\]
We shall apply Lemma \ref{lem:1} with $U=\Omega_-$, $\tilde\Gamma=OA\cup OC$ and $\phi(x,y)=\psi(\theta)$, where $\psi(\theta)$ is the solution of (\ref{ode}) described in Lemmas \ref{lem:largeb} and \ref{lem:smallb}. An easy computation shows that
\[
-\Delta \psi =\frac{c}{d^2}\psi \, .
\]
We thus obtain that
\begin{equation}
\int_{\Omega_-}|\nabla u|^2 dx\, dy \geq c\int_{\Omega_-}\frac{u^2}{d^2}dx\, dy +
\int_{\Gamma}\frac{\nabla\phi}{\phi}\cdot\vec{\nu} u^2 dS \, , \qquad u\in C^{\infty}_c(\Omega).
\label{pa}
\end{equation}
We next apply Lemma \ref{lem:1} for $\Omega_+$ and the function
$\phi_1(x,y)=d(x,y)^{\alpha}$ (we recall that $\alpha$ is the largest solution of $\alpha(1-\alpha)=c$). We note that in $\Omega_+$ the function $d(x,y)$ coincides with the distance from $AB\cup BC$ and this implies that
\[
-\Delta d^{\alpha} \geq \alpha(1-\alpha)\frac{d^{\alpha}}{d^2}\; , \qquad \mbox{ on }\Omega_+\, .
\]
(The difference of the two functions above is a positive mass concentrated on the bisector of the angle $B$).
Applying Lemma \ref{lem:1} we obtain that
\begin{eqnarray}
\int_{\Omega_+}|\nabla u|^2 dx\, dy &\geq& c\int_{\Omega_+}\frac{u^2}{d^2}dx\, dy  -\int_{\Gamma}\frac{\alpha\nabla d}{d}\cdot\vec{\nu}\, u^2 dS
\, , \qquad u\in C^{\infty}_c(\Omega).
\label{pa1}
\end{eqnarray}
Adding (\ref{pa}) and (\ref{pa1}) we conclude that
\begin{equation}
\int_{\Omega}|\nabla u|^2 dx\, dy \geq c\int_{\Omega}\frac{u^2}{d^2}dx\, dy +
\int_{\Gamma}\Big(\frac{\nabla\phi}{\phi}-\alpha\frac{\nabla d}{d}\Big)\cdot\vec{\nu}\, u^2 dS \, , \qquad u\in C^{\infty}_c(\Omega).
\label{nd}
\end{equation}
We emphasize that in the last integral the values of $\nabla\phi/\phi$ are obtained as limits from $\Omega_-$ while those of $\nabla d/d$ are obtained as limits from $\Omega_+$.

It remains to prove that the line integral in (\ref{nd}) is non-negative. For this we shall consider the different segments of $\Gamma$.

(i) The segment $\Gamma_1$ ($0\leq\theta\leq\pi/2$). Simple calculations give
\begin{equation}
\frac{\nabla\phi}{\phi}=\frac{1}{r}\frac{\psi'(\theta)}{\psi(\theta)}(-\sin\theta , \cos\theta) \, , \quad \mbox{ in } \Omega_-\, .
\label{fl1}
\end{equation}
The line $AB$ has equation $y+(x-1)\tan\gamma =0$, so $d(x,y)= (1-x)\sin\gamma  -y\cos\gamma$ on $\{P\in\Omega : d(P)=\dist(P,AB)\}$ and therefore
\begin{equation}
\alpha\frac{\nabla d}{d}=-\alpha\frac{(\sin\gamma ,\cos\gamma) }{d} \; , \qquad \mbox{ on }\Gamma_1\cup\Gamma_2 .
\label{fl2}
\end{equation}
Since $\vec{\nu}=(\sin(\gamma/2), \cos(\gamma/2))$ along $\Gamma_1$, (\ref{fl1}) and (\ref{fl2}) yield
\[
\Big(\frac{\nabla\phi}{\phi}-\alpha\frac{\nabla d}{d}\Big)\cdot\vec{\nu} =\frac{1}{r}\frac{\psi'(\theta)}{\psi(\theta)}\cos(\theta +\frac{\gamma}{2}) +
\frac{\alpha \cos(\gamma/2)}{d} \; , \qquad \mbox{ on }\Gamma_1\, .
\]
However $d(x,y)=y=r\sin\theta$ on $\Gamma_1$, so we conclude by (i) of Lemma \ref{lem:quad} (with $\omega=\gamma/2$) that
\begin{equation}
\Big(\frac{\nabla\phi}{\phi}-\alpha\frac{\nabla d}{d}\Big)\cdot\vec{\nu} =\frac{1}{r\sin\theta}\Big(g(\theta)\cos(\theta +\frac{\gamma}{2}) +
\alpha \cos(\gamma/2)\Big)\geq 0 \; , \qquad \mbox{ on }\Gamma_1\, .
\label{fl4}
\end{equation}

(ii) The segment $\Gamma_2$ ($\pi/2\leq\theta\leq\theta_0$). This is (part of) the parabola determined by the origin and the side $AB$. Applying (\ref{normal})
we obtain that the outward (with respect to $\Omega_-$) unit normal along $\Gamma_2$ is
\begin{equation}
\vec{\nu} =\frac{(\cos\theta +\sin\gamma, \sin\theta +\cos\gamma)}{\sqrt{2+2\sin(\theta+\gamma)}}.
\label{fl5}
\end{equation}
Combining (\ref{fl1}), (\ref{fl2}), (\ref{fl5}) and (ii) of Lemma \ref{lem:quad} (with $\omega=\gamma$) we obtain
\begin{equation}
\Big(\frac{\nabla\phi}{\phi}-
\alpha\frac{\nabla d}{d}\Big)\cdot\vec{\nu} =\frac{1}{r\sqrt{2+2\sin(\theta+\gamma)}}\Big(\frac{\psi'(\theta)}{\psi(\theta)}\cos(\theta+\gamma) +\alpha
[1+\sin(\theta +\gamma )]\Big)\geq 0\; , \quad \mbox{ on }\Gamma_2\, .
\label{fl6}
\end{equation}
(iii) The segment $\Gamma_3$ ($\theta_0\leq\theta\leq\beta -\pi/2$). This is (part of) the parabola determined by the origin and the side $BC$.
Now, the line $BC$ has equation
\[
 (x+T) \sin(\gamma+\delta) +y\cos(\gamma+\delta)  =0 \; ,
\]
where $(-T,0)$ is the point where the side $BC$ intersects the $x$-axis.
Applying (\ref{normal}) we thus obtain that the outward unit normal is
\[
\vec{\nu}=\frac{(\cos\theta -\sin(\gamma+\delta),\sin\theta-\cos(\gamma+\delta))}{\sqrt{2-2\sin(\theta+\gamma+\delta)}}.
\]
Hence, by (iii) of Lemma \ref{lem:quad} (with $\omega=\gamma+\delta$),
\begin{equation}
\Big(\frac{\nabla\phi}{\phi}-
\alpha\frac{\nabla d}{d}\Big)\cdot\vec{\nu} =\frac{1}{r\sqrt{2-2\sin(\theta+\gamma)}}\Big(-\frac{\psi'(\theta)}{\psi(\theta)}\cos(\theta+\gamma+\delta) +\alpha [1-\sin(\theta +\gamma+\delta )]\Big)\geq 0\; , \;\; \mbox{ on }\Gamma_3\, .
\label{fl7}
\end{equation}
(iv) The segment $\Gamma_4$ ($\beta-\pi/2\leq\theta\leq\beta$). Replacing $\theta$ by $\beta-\theta$, $\gamma$ by $2\pi -\beta-\gamma-\delta$ (the angle at $C$) and using the relation $\psi(\theta)=\psi(\beta-\theta)$, the computations become identical to those for the segment $\Gamma_1$; hence we obtain
\begin{equation}
\Big(\frac{\nabla\phi}{\phi}-\alpha\frac{\nabla d}{d}\Big)\cdot\vec{\nu}\geq 0 \; , \qquad \mbox{ on }\Gamma_4\, .
\label{fl8}
\end{equation}
The proof of the theorem is completed by combining (\ref{nd}), (\ref{fl4}), (\ref{fl6}), (\ref{fl7}) and (\ref{fl8}). $\hfill\Box$

\

{\bf\em Proof of Theorem: type A2.} In this case the curve $\Gamma$ consists of three line segments and one parabola segment. Without loss of generality we assume that starting from $\theta=0$ we first meet two line segments, then the parabola segment and then the last line segment. Then the first two line segments meet at the point $S$ with polar angle $\theta_0\leq\pi/2$ and the four components of $\Gamma$ are
\[
\Gamma_1=\{ 0\leq \theta\leq \theta_0\} , \;
\Gamma_2=\{ \theta_0\leq \theta\leq \frac{\pi}{2}\} , \;
\Gamma_3=\{ \frac{\pi}{2}\leq \theta\leq \beta-\frac{\pi}{2}\} , \;
\Gamma_4=\{\beta-\frac{\pi}{2}\leq\theta \leq\beta\}.
\]
As in the case A1, we apply Lemma \ref{lem:1} on $\Omega_-$ and $\Omega_+$ with the functions
$\phi(x,y)=\psi(\theta)$ and $\phi_1(x,y)=d(x,y)^{\alpha}$ respectively. We arrive at an inequality similar to
(\ref{nd})
and we conclude that the result will follow once we prove that
\begin{equation}
\Big(\frac{\nabla\phi}{\phi}-\alpha\frac{\nabla d}{d}\Big)\cdot\vec{\nu}\geq 0 \; , \qquad \mbox{ on }\Gamma\, .
\label{fl88}
\end{equation}

The computations along the segments $\Gamma_1$, $\Gamma_3$ and $\Gamma_4$ are identical to those for
the type A1 considered above and are omitted.

For $\Gamma_2$ we consider the point $(-T,0)$, $T>0$, where the side $BC$ intersects the $x$-axis. The distance from the line $BC$
is  $(x+T)\sin(\gamma+\delta)+y\cos(\gamma+\delta)$, therefore $\nabla d=(\sin(\gamma+\delta),\cos(\gamma+\delta))$ on $\Gamma_2$.
Moreover along $\Gamma_2$ we have $\vec{\nu}=(-\cos((\gamma+\delta)/2),\sin((\gamma+\delta)/2))$.
We also note on $\Gamma_2$ we have $d(x,y)=y=r\sin\theta$.
Combining the above we obtain that
\[
\Big(\frac{\nabla\phi}{\phi}-\alpha\frac{\nabla d}{d}\Big)\cdot\vec{\nu} =
\frac{1}{r\sin\theta}\Big[g(\theta)\sin\big(\theta +\frac{\gamma+\delta}{2}\big) +
\alpha\sin(\frac{\gamma+\delta}{2})\Big] \; , \qquad \mbox{ on }\Gamma_2,
\]
which is non-negative for $\theta\in [0,\pi/2]$ since $\gamma+\delta\leq\pi$. $\hfill\Box$

We next consider the cases where one of the two angles that are adjacent to the non-convex angle exceeds $\pi/2$. Without loss of generality we assume that $\gamma\geq\pi/2$ (the angle at the vertex $A$). We note that since $\beta_{cr} > 3\pi/2$, in this case we have $\pi\leq \beta\leq \beta_{cr}$ hence the Hardy constant is $c=1/4$.

We now divide $\Omega_+$ in two parts, $\Omega_+^A$ and $\Omega_+^C$, the parts of $\Omega_+$ with nearest boundary points on $AB$ and $BC$ respectively.
We denote by $\Gamma_*$ the common boundary of $\Omega_+^A$ and $\Omega_+^C$, that is the line segment $SB$.
We also denote by $\vec{\nu}_*$ the normal unit vector along $\Gamma_*$ which is outward with respect to $\Omega_+^A$.

\

{\bf\em Proof of Theorem: type B1.} As in the case A1, the curve $\Gamma$ is made up of four segments,
\[
\Gamma_1=\{ 0\leq \theta\leq \pi/2\} , \;
\Gamma_2=\{ \pi/2 \leq \theta\leq \theta_0 \} , \;
\Gamma_3=\{ \theta_0\leq \theta\leq \beta-\frac{\pi}{2}\} , \;
\Gamma_4=\{\beta-\frac{\pi}{2}\leq\theta \leq\beta\},
\]
where $\theta_0$ is the polar angle of the point $S$.
We use again Lemma \ref{lem:1}. On $\Omega_-$ we use the function $\phi(x,y)=\psi(\theta)$, exactly as in types A1 and A2 and we obtain that
\begin{equation}
\int_{\Omega_-}|\nabla u|^2 dx\, dy \geq \frac{1}{4}\int_{\Omega_-}\frac{u^2}{d^2}dx\, dy +
\int_{\Gamma}\frac{\nabla\phi}{\phi}\cdot\vec{\nu}u^2 dS \, , \qquad u\in C^{\infty}_c(\Omega).
\label{111}
\end{equation}
On $\Omega_+^C$ again we work as in types A1 and A2: we use the function $\phi(x,y)=d(x,y)^{1/2}$ and we obtain
\begin{equation}
\int_{\Omega_+^C}|\nabla u|^2 dx\, dy \geq \frac{1}{4}\int_{\Omega_+^R}\frac{u^2}{d^2}dx\, dy  -\frac{1}{2}\int_{\Gamma_3\cup\Gamma_4}\frac{\nabla d}{d}\cdot\vec{\nu} u^2 dS - \frac{1}{2}\int_{\Gamma_*}\frac{\nabla d}{d}\cdot\vec{\nu}_* u^2 dS \, ,
\qquad u\in C^{\infty}_c(\Omega).
\label{112}
\end{equation}
Concerning $\Omega_+^A$, we cannot use the test function $\phi=d^{1/2}$ since part (i) of Lemma \ref{lem:quad} is not valid for the full range $\pi/4<\omega<\pi/2$. So we construct a different function $\phi$. To do this we consider a second orthonormal coordinate system with cartesian coordinates denoted by $(x_1,y_1)$ and polar coordinates denoted by $(r_1,\theta_1)$. The origin $O_1$ of this system is located on the extension of the side $AB$ from $A$ and
at distance $-\cos\gamma$ from $A$, and the axes are chosen so that the point
$A$ has cartesian coordinates $(-\cos\gamma, 0)$ with respect to the new system. We note that this choice is such that
\begin{equation}
\label{113}
\mbox{the point on $\Gamma_1$ for which $\theta=\frac{\pi}{2}-\frac{\gamma}{2}$ satisfies also $\theta_1=\frac{\pi}{2}-\frac{\gamma}{2}$.}
\end{equation}

We apply Lemma \ref{lem:1} on $\Omega_+^A$ with the function $\phi_1(x,y)=\psi(\theta_1)$. This function clearly satisfies
$-\Delta \phi_1 \geq \frac{1}{4}\, d^{-2}\phi_1$, hence we obtain
\begin{equation}
\int_{\Omega_+^A}|\nabla u|^2 dx\, dy \geq \frac{1}{4} \int_{\Omega_+^A}\frac{u^2}{d^2}dx\, dy
 -\int_{\Gamma_1\cup\Gamma_2}(\frac{\nabla\phi_1}{\phi_1}\cdot\vec{\nu})u^2  \,dS
 + \int_{\Gamma_*}(\frac{\nabla\phi_1}{\phi_1}\cdot\vec{\nu}_*)u^2  \,dS
 \, \quad  u\in C^{\infty}_c(\Omega).
 \label{114}
\end{equation}
Adding (\ref{111}), (\ref{112}) and (\ref{114}) we conclude that
\begin{eqnarray}
&&\int_{\Omega}|\nabla u|^2 dx\, dy \geq \frac{1}{4}\int_{\Omega}\frac{u^2}{d^2}dx\, dy + \int_{\Gamma_1\cup\Gamma_2}\Big(\frac{\nabla\phi}{\phi}-\frac{\nabla\phi_1}{\phi_1}\Big)\cdot\vec{\nu}\, u^2 dS \nonumber\\
&&\qquad + \int_{\Gamma_3\cup\Gamma_4}\Big(\frac{\nabla\phi}{\phi}-\frac{\nabla d}{2d}\Big)\cdot\vec{\nu}\, u^2 dS +
  \int_{\Gamma_*}\Big(\frac{\nabla\phi_1}{\phi_1}-\frac{\nabla d}{2d}\Big)\cdot\vec{\nu}_*\, u^2 dS
\label{115}
\end{eqnarray}
for any $u\in C^{\infty}_c(\Omega)$. So it remains to prove that the three line integrals in (\ref{115}) are non-negative. For this we shall separately consider the different the segments $\Gamma_1$, $\Gamma_2$, $\Gamma_3$ and $\Gamma_4$ and the segment $\Gamma_*$.


(i) The segment $\Gamma_1$ ($0\leq\theta\leq\pi/2$). We have
\[
\frac{\nabla\phi}{\phi}\cdot\vec{\nu} = \frac{\psi'(\theta)}{r\psi(\theta)}\cos(\theta+\frac{\gamma}{2}) \; , \quad
\mbox{ on } \Gamma_1.
\]
and similarly
\[
\frac{\nabla\phi_1}{\phi_1}\cdot\vec{\nu} = -\frac{\psi'(\theta_1)}{r_1 \psi(\theta_1)}\cos(\theta_1-\frac{\gamma}{2}) \; , \quad
\mbox{ on }  \Gamma_1.
\]
However we have $r_1\sin\theta_1 =r\sin\theta$ along $\Gamma_1$, so recalling definition (\ref{g}) we see that it is enough to prove the inequality
\begin{equation}
g(\theta)\cos(\theta+\frac{\gamma}{2}) +g(\theta_1)\cos(\theta_1-\frac{\gamma}{2})\geq 0 \; , \qquad \mbox{ on }\Gamma_1\, .
\label{222}
\end{equation}
Recalling (\ref{113}) and  applying the sine law we obtain that along $\Gamma_1$ the polar angles $\theta$ and $\theta_1$ are related by
\begin{equation}
\cot \theta_1 =-\cos\gamma \cot\theta +\sin\gamma \; .
\label{223}
\end{equation}
{\bf Claim.} There holds
\begin{equation}
\label{claim1}
 \theta_1\geq  \theta +\gamma -\pi \; , \qquad \mbox{ on }\Gamma_1 \, .
\end{equation}
{\em Proof of Claim.} We fix $\theta\in [0,\pi/2]$ and
the corresponding $\theta_1=\theta_1(\theta)$. If $\theta+\gamma -\pi\leq 0$, then (\ref{claim1}) is obviously true, so we assume that $\theta+\gamma -\pi\geq 0$. Since $0\leq\theta+\gamma -\pi\leq \pi/2$ and $0\leq\theta_1\leq\pi/2$, (\ref{claim1}) is written equivalently $\cot\theta_1\leq\cot(\theta +\gamma -\pi)$; thus, recalling (\ref{223}), we conclude that to prove the claim it is enough to show that
\[
 -\cos\gamma\cot\theta +\sin\gamma\leq\cot(\theta+\gamma) \; , \quad
\pi-\gamma\leq\theta\leq\frac{\pi}{2},
\]
or, equivalently (since $\pi\leq\theta+\gamma\leq 3\pi/2$),
\begin{equation}
-\cos\gamma\cot^2\theta +(-\cos\gamma\cot\gamma -\cot\gamma+\sin\gamma)\cot\theta + 1+\cos\gamma \geq 0\; , \qquad
\pi-\gamma\leq\theta\leq\frac{\pi}{2}.
\label{eq:new}
\end{equation}
The left-hand side of (\ref{eq:new}) is an increasing function of $\cot\theta$ and therefore takes its least value at $\cot\theta =0$.
Hence the claim is proved.

For $0\leq\theta\leq\pi/2-\gamma/2$ (\ref{222}) is true since all terms in the left-hand side are non-negative. So let
$\pi/2-\gamma/2\leq\theta\leq\pi/2$ and $\theta_1=\theta_1(\theta)$. From (\ref{223}) we find that
\begin{eqnarray*}
\frac{d\theta_1}{d\theta}-1&=&-\frac{\cos\gamma (1+\cot^2\theta)+1+\cot^2\theta_1}{1+\cot^2\theta_1} \\
&=& -\frac{ 1+\sin^2\gamma +\cos\gamma -2\sin\gamma\cos\gamma\cot\theta +\cos\gamma(1+\cos\gamma)\cot^2\theta}{1+\cot^2\theta_1}.
\end{eqnarray*}
The function
\[
h(x):= 1+\sin^2\gamma +\cos\gamma -2\sin\gamma\cos\gamma x +\cos\gamma(1+\cos\gamma)x^2
\]
is a concave function of $x$. We will establish the positivity of $h(\cot\theta)$ for $\pi/2-\gamma/2
\leq \theta\leq\pi/2$. For this it is enough to establish the positivity at the endpoints. At $\theta=\pi/2$ positivity is obvious, whereas
\[
h(\tan(\frac{\gamma}{2}))=1+\sin^2\gamma+\cos\gamma -2\cos\gamma\sin^2\frac{\gamma}{2}\geq 0.
\]
From (\ref{113}) we conclude that $\theta_1\leq\theta$ for $\pi/2-\gamma/2\leq\theta\leq\pi/2$.

We next apply Lemma \ref{lem:g}. We obtain that for $\pi/2-\gamma/2\leq\theta\leq\pi/2$,
\begin{eqnarray*}
g(\theta)\cos(\theta+\frac{\gamma}{2}) +g(\theta_1)\cos(\theta_1-\frac{\gamma}{2})
&\geq& g(\theta) [\cos(\theta+\frac{\gamma}{2}) +\cos(\theta_1-\frac{\gamma}{2})]\\
&=&2g(\theta)\cos(\frac{\theta+\theta_1}{2}) \cos (\frac{\theta-\theta_1+\gamma}{2})\\
&\geq& 0,
\end{eqnarray*}
where for the last inequality we made use of the claim. Hence (\ref{222}) has been proved.

(ii) The segment $\Gamma_2$ ($\frac{\pi}{2}\leq\theta\leq\theta_0$). Computations similar to those that led to (\ref{fl6}) together with the fact that $r=r_1\sin\theta_1$ on $\Gamma_2$ give that along $\Gamma_2$ we have
\begin{eqnarray}
&&\Big(\frac{\nabla\phi}{\phi}-\frac{\nabla\phi_1}{\phi_1}\Big)\cdot\vec{\nu} \nonumber\\
&=&\frac{1}{\sqrt{2+2\sin(\theta+\gamma)}}\Big[ \frac{f(\theta)}{r} \cos(\theta+\gamma)
- \frac{f(\theta_1)}{r_1}[\sin(\theta_1-\theta -\gamma)-\cos\theta_1] \Big] \label{night}\\
&=&\frac{1}{r\sqrt{2+2\sin(\theta+\gamma)}}\Big[ f(\theta)\cos(\theta+\gamma)- f(\theta_1)\sin\theta_1[\sin(\theta_1-\theta -\gamma)-\cos\theta_1] \Big].\nonumber
\end{eqnarray}
Now, simple geometry shows that along $\Gamma_2$ the angles $\theta$ and $\theta_1$ are related by
\begin{equation}
\cot\theta_1 =-\cos(\theta+\gamma).
\label{224}
\end{equation}
It follows that
\[
\sin\theta_1 [\sin(\theta_1-\theta-\gamma)-\cos\theta_1]= \frac{\cos(\theta+\gamma) [2+\sin(\theta+\gamma)]}{1+\cos^2(\theta+\gamma)},
\;\;\;\;\mbox{ along $\Gamma_2$}\, .
\]
Since $\cos(\theta+\gamma)\leq 0$, (\ref{224}) and Lemma \ref{lem:sss}
imply that $(\nabla\phi/\phi -\nabla\phi_1/\phi_1)\cdot\vec{\nu}\geq 0$ along $\Gamma_2$, as required.

(iii) The segments $\Gamma_3$ and $\Gamma_4$ ($\theta_0\leq\theta \leq \beta$). Since $\zeta <\pi/2$, the change $\theta \leftrightarrow
\beta -\theta$ reduces this case to that of the segments $\Gamma_2$ and $\Gamma_1$ respectively for a quadrilateral of type A1, already considered above.

(iv) The segment $\Gamma_*$. The contribution from $\Omega_+^A$ is
\[
\frac{\nabla\phi_1}{\phi_1}\cdot \vec{\nu}_* = \frac{f(\theta_1)}{r_1}\cos(\theta_1+\frac{\delta}{2})\geq 0\, , \qquad 
\mbox{ on }\Gamma_*,
\]
since $\theta_1\leq \gamma/2$, by construction of the new coordinate system and $\gamma+\delta<\pi$. Given that the contribution from $\Omega_+^C$ is positive, the proof is complete.

\

{\bf\em Proof of Theorem: type B2.} As in the case of type A2, there exists an angle $\theta_0\leq\pi/2$ such that the four segments of $\Gamma$ are
\[
\Gamma_1=\{ 0\leq \theta\leq \theta_0\} , \;
\Gamma_2=\{ \theta_0\leq \theta\leq \frac{\pi}{2}\} , \;
\Gamma_3=\{ \frac{\pi}{2}\leq \theta\leq \beta-\frac{\pi}{2}\} , \;
\Gamma_4=\{\beta-\frac{\pi}{2}\leq\theta \leq\beta\}.
\]
So $\Gamma_3$ is a parabola segment while $\Gamma_1$, $\Gamma_2$ and $\Gamma_4$ are line segments.
We define the sets $\Omega_+^A$, $\Omega_+^C$ and the vector $\vec{\nu}_*$ as in the case of type B1 and apply Lemma \ref{lem:1} with the same functions, that is
$\psi(\theta)$ on $\Omega_-$, $d(x,y)^{1/2}$ on $\Omega_+^C$ and $\psi(\theta_1)$ on $\Omega_+^A$ (where we use exactly the some construction for the coordinate system $(x_1,y_1)$).

The computations along $\Gamma_1$, $\Gamma_3$ and $\Gamma_4$ are identical to those for the type B1 and are omitted. On $\Gamma_2$ we have, as in the case of subtype A2,
\[
\Big(\frac{\nabla\phi}{\phi}-\alpha\frac{\nabla d}{d}\Big)\cdot\vec{\nu} =
\frac{1}{r\sin\theta}\Big[g(\theta)\sin\big(\theta +\frac{\gamma+\delta}{2}\big) +\half\sin(\frac{\gamma+\delta}{2})\Big] \geq 0\, ,
\]
since $\gamma+\delta\leq\pi$. Finally, the computations along $\Gamma_*$ are identical to the corresponding computations for the case $B1$. This completes the proof.

\

{\bf\em Proof of Theorem: Type B3.} In this case there exist angles $\theta_0,\theta_0'$ with
\[
\frac{\pi}{2}\leq \theta_0  <\theta_0' \leq \beta -\frac{\pi}{2}
\]
such that the four segments of $\Gamma$ are
\[
\Gamma_1=\{ 0\leq \theta\leq \frac{\pi}{2}\} , \;
\Gamma_2=\{ \frac{\pi}{2}\leq \theta\leq \theta_0\} , \;
\Gamma_3=\{ \theta_0\leq \theta\leq  \theta_0' \} , \;
\Gamma_4=\{\theta_0'\leq\theta \leq\beta\}.
\]
So $\Gamma_2$ is a parabola segment while $\Gamma_1$, $\Gamma_3$ and $\Gamma_4$ are line segments. To proceed,
we define the sets $\Omega_+^A$, $\Omega_+^C$ and the vector $\vec{\nu}_*$ as in the cases B1 and B2 and apply Lemma \ref{lem:1} with the same functions, that is $\psi(\theta)$ on $\Omega_-$, $d(x,y)^{1/2}$ on $\Omega_+^C$ and $\psi(\theta_1)$ on $\Omega_+^A$, where again we use exactly the some construction for the coordinate system $(x_1,y_1)$.

The computations for the line segments $\Gamma_1$ and $\Gamma_4$ and for the parabola segment $\Gamma_2$ are identical to those for a quadrilateral of type B1 and are omitted. We next consider the line segment $\Gamma_3$ whose points are equidistant from the sides $AB$ and $OC$. Calculations similar to those above give that
\[
\Big(\frac{\nabla\phi}{\phi}-\frac{\nabla \phi_1}{\phi_1}\Big)\cdot\vec{\nu} =
\frac{1}{r\sin\theta}\Big[g(\theta)\sin\big(\frac{\beta-\gamma}{2}-\theta\big) +
g(\theta_1)\sin(\frac{\beta+\gamma}{2}-\theta_1)\Big] \; , \;\;\mbox{ on }\Gamma_3.
\]
Now, it follows by construction that
\[
\theta \geq \frac{\pi}{2} \geq \frac{\beta+\gamma-\pi}{2} \geq\theta_1 \;\; , \qquad \mbox{ on }\Gamma_3.
\]
Since $0 < (\beta+\gamma)/2  -\theta_1  <\pi$, by the monotonicity of $g$ we have
\begin{eqnarray*}
\Big(\frac{\nabla\phi}{\phi}-\frac{\nabla \phi_1}{\phi_1}\Big)\cdot\vec{\nu} &\geq& \frac{g(\theta)}{r\sin\theta}
\Big[\sin\big(\frac{\beta-\gamma}{2}-\theta\big) + \sin(\frac{\beta+\gamma}{2}-\theta_1)\Big] \\
&=& \frac{2g(\theta)}{r\sin\theta}\sin\big( \frac{\beta-\theta-\theta_1}{2}\big)\cos\big( \frac{\gamma+\theta -\theta_1}{2}\big).
\end{eqnarray*}
Since $0< \beta -\theta-\theta_1  <2\pi$, the last sine is positive. It is also clear that $\gamma +\theta-\theta_1 >0$.
Hence the proof will be complete if we establish the following

{\bf Claim:} There holds
\begin{equation}
\theta_1 \geq \theta+\gamma-\pi \; \; , \qquad \mbox{ on }\Gamma_3.
\label{abcd}
\end{equation}
{\em Proof of Claim.}  Simple geometry shows that along $\Gamma_3$ the polar angles $\theta$ and $\theta_1$ are related by
\[
\cot\theta_1 =-\cos(\beta+\gamma) \cot(\beta-\theta) -\sin(\beta+\gamma) \; .
\]
and  $[\theta_0,\theta_0']\subset [\pi/2, \beta - \pi/2]\subset  [\pi/2, (\beta -\gamma+ \pi)/2] $. We will actually establish (\ref{abcd}) for the larger range $\pi/2\leq\theta \leq (\beta -\gamma+ \pi)/2$.

For this, we initially observe that for $\theta=(\beta-\gamma+\pi)/2$ inequality (\ref{abcd}) holds as an equality. Therefore the claim will be proved if we establish that
\[
\frac{d\theta_1}{d\theta} -1 \leq 0 \; , \quad\quad \frac{\pi}{2}\leq \theta\leq \frac{\beta-\gamma+\pi}{2}.
\]
However, we easily come up to
\begin{eqnarray*}
\frac{d\theta_1}{d\theta} -1 &=& -\frac{  \cos(\beta+\gamma)(\cos(\beta+\gamma) -1)\cot^2(\beta-\theta) +2\sin(\beta+\gamma)\cos(\beta+\gamma)\cot(\beta-\theta)}{1+\cot^2\theta_1} \\
&& -\frac{1 +\sin^2(\beta+\gamma)-\cos(\beta+\gamma)}{1+\cot^2\theta_1}.
\end{eqnarray*}
The function
\[
h(x):=  \cos(\beta+\gamma)(\cos(\beta+\gamma) -1)x^2 +2\sin(\beta+\gamma)\cos(\beta+\gamma)x +1 +\sin^2(\beta+\gamma)-\cos(\beta+\gamma)
\]
is a concave function of $x$. We will establish the positivity of $h(\cot(\beta-\theta))$, $\pi/2\leq \theta\leq (\beta-\gamma+\pi)/2$, and for this it is enough to establish positivity at the endpoints. A simple computation shows that
\[
h( \cot(\beta- \frac{\beta -\gamma+\pi}{2})) =2\tan^2(\frac{\beta+\gamma}{2}).
\]
At the other endpoint we have
\begin{eqnarray*}
 h(\cot(\beta-\frac{\pi}{2}))&=& \cos(\beta+\gamma)(\cos(\beta+\gamma) -1)\tan^2\beta -\\
&& -2\sin(\beta+\gamma)\cos(\beta+\gamma)\tan\beta+1+\sin^2(\beta+\gamma)-\cos(\beta+\gamma)\\
& =&\frac{2 \sin^2(\frac{\beta+\gamma}{2})}{\cos^2\beta}\big[  1+ \cos(2\beta)\cos^2(\frac{\beta+\gamma}{2})  \big] \\
&& - \frac{\sin(\beta+\gamma)}{2\cos^2\beta}\big(  \sin(\beta-\gamma) +\sin(2\beta) \cos(\beta+\gamma)\big)\\
&\geq& 0,
\end{eqnarray*}
since $3\pi/2 \leq \beta+\gamma\leq 2\pi$ and $0\leq \beta-\gamma\leq \pi$. Hence the claim is proved and therefore the total contribution along $\Gamma_3$ is non-negative.

It finally remains to establish that the total contribution along $\Gamma_*$ is non-negative. As in type B1 the contribution from $\Omega_+^A$ is
\[
\frac{\nabla\phi_1}{\phi_1}\cdot \vec{\nu}_* = \frac{f(\theta_1)}{r_1}\cos(\theta_1+\frac{\delta}{2}).
\]
This is is non-negative since $\theta_1 <(\beta+\gamma-\pi)/2$ and $\beta+\gamma+\delta <2\pi$.
This completes the proof.
$\hfill\Box$





\end{document}